\def\E{\end{document}}
\documentclass[11pt]{article}
\usepackage{amssymb,amsmath}

\begin{document}
\title{}
\date{}
\author{}
\title{\bf Feedback controllability for  blowup points  of  heat equation
\footnote{This work was partially supported by the National Natural
Science Foundation of China under grant 12071067, and
 National Key R\&D Program of China under grant 2020YFA0714102.}}

\date{}
\author{Ping Lin\footnote{  School of
Mathematics \& Statistics, Northeast Normal University, Changchun
130024,  P. R. China. E-mail address:
linp258@nenu.edu.cn.} \qquad Hatem Zaag\footnote{
Universit\'e Sorbonne Paris Nord,
LAGA, CNRS (UMR 7539), F-93430, Villetaneuse, France. Email:
hatem.zaag@univ-paris13.fr.
}
\\
\\
 }
\maketitle
\newtheorem{theorem}{Theorem}[section]
\newtheorem{definition}{Definition}[section]
\newtheorem{lemma}{Lemma}[section]
\newtheorem{proposition}{Proposition}[section]
\newtheorem{corollary}{Corollary}[section]
\newtheorem{remark}{Remark}[section]
\renewcommand{\theequation}{\thesection.\arabic{equation}}
\catcode`@=11 \@addtoreset{equation}{section} \catcode`@=12
{{\bf Abstract.}
This paper concerns a  controllability problem for blowup points on heat equation.  It can be described as follows: In the absence of control, the solution to the linear
heat system
 globally exists in a bounded domain $\Omega$.  While,  for a given time $T>0$ and a point $a$ in this domain, we find a feedback control, which is acted on an internal
subset $\omega$ of this domain,
such that the corresponding solution to this system blows up at time $T$ and holds unique point $a$.  We  show   that  $a\in \omega$ can be the unique blowup point of
the corresponding solution  with  a certain feedback control, and  for any feedback control, $a\in \Omega\setminus \overline{\omega}$ could not  be the unique blowup point.\\



\section{Introduction}

 In nature and in practical applications, blowup phenomenon exists widely. This phenomenon can be described by some
nonlinear evolution differential  equations, which
 have been studied in many references (see, for instance, \cite{Abdelhedi}, \cite{Bricmont}, \cite{Filippas}-\cite{Guo-Hu}, \cite{Herrero-1}-\cite{Khenissy}, \cite{Zaag2}-\cite{Nguyen1},
 \cite{J.J.L.}-\cite{chunpeng3}).
 Roughly speaking, blowup is a conception  which means that a solution is
unbounded in finite time. In certain cases,
the blowup of a solution is desired.
For instance,
the dramatic increase in temperature  leads to the ignition of a chemical
reaction.
 However,  solutions to linear partial differential  equations  without control
generally globally exist. It is naturally interesting to find  feedback controls
 to these equations such that the corresponding solutions blow up in finite time and at given place.

This paper concerns a  controllability problem for blowup points on heat equations. It can be described as follows: In the absence of control, the solution to the linear
heat system
 globally exists in a bounded domain. While,  for a given time $T>0$ and a point $a$ in this domain, we find a feedback control, which is acted on an internal
subset of this domain,
such that the corresponding solution to this system blows up at time $T$ and holds unique point $a$. More precisely, we can describe our problem as follows.

Let $\Omega$ be a bounded domain of $\mathbb{R}^n$ with smooth boundary $\partial \Omega$. In this paper,
we
consider the following   control system,
\begin{eqnarray}\label{e1.1}
\left\{\begin{array}{ll} y_t-
\Delta y=\chi_\omega u,&x\in \Omega,\ t>0,\\
y=0,& x\in \partial\Omega,\  t>0,\\
\displaystyle  y(x,0)=y_0(x),&
x\in \Omega.
\end{array}\right.
\end{eqnarray}
Here the control $u$ acts on a nonempty and open subset
$\omega\subset\Omega$, and  $\chi_\omega$ is the characteristic
function of the set $\omega$.

It is well known that if the control $u\equiv 0$, then for any $y_0\in H_0^1(\Omega)$, the corresponding solution
to (\ref{e1.1}) globally exists.

The following definition refers to a blowup point of a solution  to  (\ref{e1.1}) with a feedback  control.
\begin{definition}
Let $T>0$ and $a\in \Omega$. We say that $T$ is the blowup time and $a$ is  a blowup point of the  solution $y$ to system (\ref{e1.1}) with a a feedback control
 $u$,  if there
  are  sequence $\{a_j\}_{j=1}^{\infty}$ with $a_j\to a$ and  sequence $\{t_j\}_{j=1}^{\infty}$ with $t_j \to T$ such that $|y(a_j,t_j)|\to +\infty$, as $j\rightarrow\infty$.
\end{definition}


The problem we studied in this paper is as follows.\\ \\
\textbf{Problem ($P$)}  Given an initial data $y_0\in H_0^1(\Omega)$,  $a\in \Omega$ and $T>0$,    we find a  feedback control $u$
 such that $T$ is the blowup time of the corresponding solution $y$ to (\ref{e1.1}), and $y$  has unique blowup point  $a$.

Our problem differs with classical controllability ones (see, for instance, \cite{zuazuax1}, \cite{CORON3}-\cite{Fernandez-Cara}, \cite{zuazuax3},
\cite{zuazuax4}, \cite{chunpeng4}, \cite{chunpeng5}). The target of our controllability problem is ``infinity",  which is outside the state spaces
of the solutions; While the targets of classical controllability problems
 are within the state spaces.
 On the other hand, feedback control can form a closed-loop system and play an effective role of control. Its characteristics are: to make timely response to the objective effect caused by each
 step of the implementation process of the plan decision, and accordingly adjust and modify the next step of the implementation plan, so that the implementation of the plan decision
 and the original plan itself can achieve dynamic coordination.

So far, there are few papers on the controllability of equations with the property of blowup.
Zuazua et al. \cite{Doubova} and \cite{Fernandez-Cara} considered the controllability of weakly blowing up semilinear parabolic
equations with open-loop controls. They showed that the systems considered are null and approximately controllable at any time. In these references, blowup  occurs in the absence of control; While
the solution can be steered to be zero or to be sufficiently approximate to a given target in the state space at given time
 by using controls.
The aim  in \cite{Doubova, Fernandez-Cara} is to prevent blowup
by controls, which is totally different with the intention of making the solution blow up in this paper.

As for feedback blowup controllability,  Lin \cite{Lin} considered the blowup   controllability of   heat equation with feedback controls.  It was proved in \cite{Lin}
that for any initial data in $H^{1}_{0}(\Omega)$  and for any time $T>0,$ there exist a number $p$ with $1<p<\infty$ and a feedback control acting on an internal
subset of the space domain such that the $L^{p+1}$ norm of the corresponding solution for (\ref{e1.1}) blows up at $T.$ Lin \cite{Lin2} et. al derived a global exact blowup controllability
for ordinary differential system $y'(t)=Ay(t)+Bu(t)$ in the case that $(A,B)$ is null controllable,
 $A$ and $B$ are time-invariant matrix. More precisely, for any initial data in $\mathbb{R}^n$  and for any time $T>0,$  one can find a feedback control $u$ to make
 the solution $y(\cdot)$ to this ODE system  blow up at time $T>0,$  i.e., $\lim\limits_{t\rightarrow T}|y(t)|_{\mathbb{R}^n}=+\infty.$

Sometimes people need more accurate blowup. For instance, in the mining process, people expect that blowup happens in a specified time and at a given place. It need more explicit theory than those obtained in \cite{Lin}
and \cite{Lin2}. The problem of blowup point controllability  studied in this paper is to meet this need, which contains  more complex theoretical analysis and techniques.

Take $H=L^2(\Omega)$,  $D(A)=H_0^1(\Omega)\cap H^2(\Omega)$, $A=\Delta: D(A)\subset H\rightarrow H$, and
$B: H\rightarrow H$, $Bu=\chi_\omega u$.  Denote by $\Sigma^+(H)$ the Banach space of all symmetric and positive operators acting in $H$.
The following theorem is our main result of  blowup controllability with a feedback control. It establishes the existence of  a blowup solution  with   prescribed profile in given time and at given point.

\begin{theorem} \label{Main Tho}     For any $y_0\in H_0^1(\Omega)$,  any $a\in \omega$ and any $T>0$, there exist $T_1\in (0, T/2)$ and $\widetilde{y}_0\in C_0^\infty(\omega)$ such that   the  solution
$y$ to  (\ref{e1.1}) with the following feedback control $$u(x,t):=\left\{\begin{array}{ll}-B^*P(t)({y}(t)-\widetilde{y}_0\big)(x)-\Delta \widetilde{y}_0(x),\ &(x,t)\in \Omega\times  (0,T-T_1), \\
|y|^{p-1}y(x,t), \ &(x,t)\in \Omega\times [T-T_1,T),
\end{array}\right.$$
exists on $[0,T)$, $T$ is the blowup time of $y$ and $y$
has unique blowup point $a$, where $p$ could be any real number with $p>1$.
Here, $P\in C_S\big([0,T-T_1); \Sigma^+(H)\big)$ is the unique mild solution to the following Riccati system,
 \begin{eqnarray}\label{13-2}
\left\{\begin{array}{ll} P'(t)+A^*P(t)+P(t)A-P(t)BB^*P(t)=0 \ \mbox{on}\ [0,T-T_1)
,\\
\lim\limits_{(s,z)\rightarrow(T-T_1,z_0)}\langle P(s)z, z \rangle=+\infty,\ \mbox{for each} \ z_0\in H\
\mbox{and}\ z_0\neq 0,
\end{array}\right.
\end{eqnarray}
and $P(\cdot)$ holds the property that   $\lim\limits_{t\rightarrow T-T_1}\langle P(t)z(t), z(t) \rangle=0$ for every mild solution $z$ of the
state system $z'=Az+Bv$, $z(t_0)=z_0$  with $0\leq t_0<T-T_1$, $z(T-T_1)=0$ and $v\in L^2(t_0,T-T_1;H)$.

Moreover, for all $R>0$,
 \begin{eqnarray}\nonumber
\sup\limits_{\big\{|x-a|\leq R
\sqrt{(T-t)|\log(T-t)|}\big\}}\Big|(T-t)^{\frac{1}{p-1}}y(x,t)-f\Big(\frac{x-a}{\sqrt{(T-t)|\log(T-t)|}}\Big)\Big|\rightarrow 0\ \mbox{as} \ t\rightarrow T,
\end{eqnarray}
where
\begin{eqnarray}\label{f}
f(\eta)=\Big(p-1+\frac{(p-1)^2}{4p}|\eta|^2\Big)^{-\frac{1}{p-1}}, \ \forall\ \eta\in \mathbb{R}.
\end{eqnarray}
\end{theorem}

\begin{remark}\label{rem1} We denote by $C_S([0,T-T_1);\Sigma^+(H))$ the set of all mappings $S:[0,T-T_1)\rightarrow \Sigma^+(H)$ such that $S(\cdot
)z_0$ is continuous on $[0,T-T_1)$ for each $z_0\in H$.  $P\in C_S([0,T-T_1); \Sigma^+(H))$ is called a mild solution to system
(\ref{13-2}) if for each $\delta\in (0, T-T_1)$, $P$ satisfies
 \begin{eqnarray}&&\nonumber P(T-T_1-\delta-t)z_0=e^{tA^*}P(T-T_1-\delta)e^{tA}z_0\\&&-\int_0^{t}e^{(t-s)A^*}P(T-T_1-\delta-s)BB^*P(T-T_1-\delta-s)e^{(t-s)A}z_0ds,\nonumber\end{eqnarray}
for each $t\in[0,T-T_1-\delta]$ and $z_0\in H$, and the second equality of (\ref{13-2}) holds.
\end{remark}

 \begin{theorem} \label{Tho 2}  Suppose that $y_0\in H_0^1(\Omega)\bigcap L^\infty(\Omega)$.  If $y$ is a corresponding solution to system (\ref{e1.1}) for some feedback control and belongs to   $C([0,t_{max});L^\infty(\Omega))$, then any  $a\in \Omega\setminus \overline{\omega}$ could not  be the unique blowup point of $y$. Here, $[0,t_{max})$ denotes the maximal interval of existence of $y$.
\end{theorem}
\begin{remark}\label{rem2} By Theorem \ref{Main Tho} and Theorem \ref{Tho 2}, we have that  $a\in \omega$ can be the unique blowup point of
the solution to system (\ref{e1.1}) with  a certain feedback control, while  for any feedback control, $a\in \Omega\setminus \overline{\omega}$ could not  be the unique blowup point of the solutions to system (\ref{e1.1}).  However, the case $a\in\partial \omega$ is open now. The methods we used to prove Theorem \ref{Main Tho} and Theorem \ref{Tho 2} are not applicable to this case. We leave it to be studied in the future. \end{remark}

The techniques to prove our main result Theorem \ref{Main Tho}  can be described as the following two steps.

 First, for any $T>0$ and any $a\in \omega$, we will find a special initial data and a feedback control such
that the corresponding  solution blows up in time $T_1$ with $0<T_1<T/2$ at one blowup point $a$, and with the prescribed blowup profile.

We will draw on the methods  that many authors studied the blowup profile of the equation $y_t-\Delta y=|y|^{p-1}y$ with $p>1$ ($*$)  in the first step. The description of the asymptotic blowup behavior, locally near a given blowup point is a main direction for this equation (see, for instance, \cite{Abdelhedi}, \cite{Bricmont},  \cite{Herrero-1}-\cite{Khenissy}, \cite{Zaag2}-\cite{Nguyen1},
 \cite{J.J.L.}, \cite{J.J.L.-1}).
   Stimulated by the blowup rate estimate,  the notion of blowup profile was introduced by Herrero and
Vel\'{a}zquez  \cite{Herrero-1, Herrero}, Vel\'{a}zquez  \cite{J.J.L., J.J.L.-1}, Filippas and Kohn \cite{Filippas} and Filippas and Liu
 \cite{Filippas-1}. The
selfsimilar change of variables is particularly well adapted to the study of the blow-up profile.

 The first step uses the ideas developed by Bricmont and Kupianinen \cite{Bricmont} and  Merle and Zaag \cite{Zaag1} to construct a blowup solution for the semilinear heat equation ($*$), and Mahmoudi, Nouaili and Zaag \cite{Zaag2} to
construct a periodic solution to ($*$) in one space dimensional, which  blows up in time $T$  at one blowup point $a$.

More precisely, in the first step (see Section 2),
we will take the feedback control $u=\chi_\omega |y|^{p-1}y$ ($p>1$) and make use of the similar techniques in Mahmoudi, Nouaili and  Zaag \cite{Zaag2}.
Indeed,  the equation $y_t-\Delta y=\chi_\omega|y|^{p-1}y$ with ($p>1$) ($**$) will be considered in Section 2. We will  prove that for any  $a\in \omega$
and any $T>0$, one can find a special initial data  such
that the corresponding  solution   blows up in time $T_1$ with $0<T_1<T/2$ at one blowup point $a$, and with the prescribed blowup profile. Our work will be divided into two parts: the blowup region and the regular region.

$\bullet$ In the blowup region, we reduce the question to a finite-dimensional problem. Similarity variables will be used to control the solution near the profile.

$\bullet$ In the regular region, we directly use the standard parabolic estimates.

We proceed by contradiction to solve the finite-dimensional problem and complete the proof of this step by index theory.

The main difference between  the problem studied in this step and that in Mahmoudi, Nouaili and  Zaag \cite{Zaag2} is that the nonlinear term of ($**$) is supported in $\omega$. Fortunately, since $a\in \omega$,  suitable cut-off functions can make us to reasonably divide the blowup region and the regular region as in  \cite{Zaag2} and to choose  an initial data $y_0\in C_0^\infty(\Omega)$ of support in $\omega$. This leads us to get proper estimates in  both regions.

 Then, in the second step, just as we did in Lin \cite{Lin}, we will prove that for each initial data, using the feedback null controllability results for linear heat equations obtained in \cite{Sirbu}, there exists a feedback
control such that the corresponding solution can  reach the above-mentioned special initial data constructed in the first step at $T-T_1$.  Combining these two steps,
we can get our desired global blowup controllability for blowup points with feedback controls.

 The rest of this paper is structured as follows. In Section 2, we will construct a blowup solution for ($**$) with   prescribed profile.
Section 3 and Section 4  will give the proof of Theorem 1.1 and Theorem 1.2, respectively.
 In Appendix, we will give the proof of a preliminary lemma.
\section{Construct a blowup solution  with   prescribed profile}
Let $p>1$ be arbitrary but fixed. We consider
the following system,
\begin{align}\label{1.6}
\left\{\begin{array}{ll} y_t-\Delta
y=\chi_\omega|y|^{p-1}y, \ \ \ \ \ \ x\in \Omega,\ t>0,\\
y=0,\ \ \ \ \ \ \ \ \ \ \ \ \ \ \ \ \ \ \
\ \ \ \ \ \ \ x\in \partial\Omega, \  t>0,\\
\displaystyle  y(x,0)=y_0(x),\ \ \ \ \ \ \ \ \ \ \ \ \ \  x\in
\Omega.
\end{array}\right.
\end{align}

In this section, we will  construct a blowup solution for (\ref{1.6}) with   prescribed profile.
 This will play an important role in the proof of our main result. More precisely, we have the following theorem.

\begin{theorem}\label{Pros1} For any $a\in \omega$, there exists $T_0>0$ such that for any $T\in (0,T_0)$, there exists an initial data $y_0\in C_0^\infty(\Omega)$ of support in $\omega$ such that the corresponding   solution
$y$  to  (\ref{1.6}) exists on $[0,T)$, $T$ is the blowup time of $y$ and $y$
has unique blowup point $a$. Moreover, for all $R>0$,
 \begin{eqnarray}\label{1.3}
\sup\limits_{\big\{|x-a|\leq R
\sqrt{(T-t)|\log(T-t)|}\big\}}\Big|(T-t)^{\frac{1}{p-1}}y(x,t)-f\Big(\frac{x-a}{\sqrt{(T-t)|\log(T-t)|}}\Big)\Big|\rightarrow 0,
\end{eqnarray}
as $t\rightarrow T$, where
\begin{eqnarray}\label{f}
f(\eta)=\Big(p-1+\frac{(p-1)^2}{4p}|\eta|^2\Big)^{-\frac{1}{p-1}}, \ \forall\ \eta\in \mathbb{R}.
\end{eqnarray}
\end{theorem}

The following four subsections will be devoted to the proof of Theorem {\ref{Pros1}}.  The techniques   are stimulated by \cite{Zaag1} and \cite{Zaag2}  for constructing  stable solutions of the type  $y_t-\Delta
y=|y|^{p-1}y$ with prescribed
profiles. In the following, we will give the proof when $n=1$, for simplicity, and that the proof for $n\geq2$ is the same, with small adaptations, as one can see from  the paper \cite{Nguyen2} done for the standard heat equation.
In subsection 2.1, we will give a formulation of our problem. In subsection 2.2 and subsection 2.3, we will reduce an infinite dimensional problem into finite
 dimensional one. In subsection 2.4, we will complete the proof of Theorem \ref{Pros1} by topological argument.

\subsection{Formulation of the problem}

Suppose that the initial data $y_0\in C_0^\infty(\Omega)$, and  the solution $y$ to (\ref{1.6}) exists on $[0,T)$. Then it is well known that for any $\delta\in(0,T)$,  $y$ is in the space $L^\infty(\Omega\times(0,T-\delta))$. Then, we have $|y|^{p-1}y(\cdot)\in L^\infty(\Omega\times(0,T-\delta))$.
Thus, by the standard $L^p$ theory of linear parabolic equations, it holds that $y\in W_q^{2,1}(\Omega\times(0,T-\delta))$ for any $q$ with $1<q<\infty$.
By  Sobolev embedding theorem,  the solution $y$ is in $C^{\alpha,\alpha/2}(\overline{\Omega\times(0,T-\delta)})$ for some $\alpha$ depending
on $n$ with $0<\alpha<1$.

Let $a\in \omega$ and $T>0$.  Let $\chi_0\in C_0^\infty(\mathbb{R},[0,1])$ with
\begin{align}\label{1.7}\chi_0(\xi)=\left\{\begin{array}{ll} 1, \ |\xi|\leq 1,\\0, \ |\xi|\geq 2.
\end{array}\right.\end{align}
Our analysis can be divided into the following two parts: the blowup region and regular region.

In the regular region $y$, we define $\overline{y}$ by
\begin{align}\overline{y}(x,t)=& y(x,t)\overline{\chi}(x),\ x\in \Omega,\ t\geq 0,\end{align}
where for any $\xi\in \Omega$,  $\overline{\chi}(\xi)=1-\chi_0(\frac{4(\xi-a)}{\varepsilon_0})$,
with $\varepsilon_0>0$ will be fixed sufficiently small  later.
Then, $\overline{y}$ satisfies the following equation:
\begin{align}\label{1.8}
\partial_t \overline{y}=\partial_{xx}\overline{y}+\chi_\omega|y|^{p-1}\overline{y}-2\overline{\chi}'\partial_x y-\overline{\chi}{'}{'}y.
\end{align}

In the blowup region of $y$, we make the following similarity transformation of system (\ref{1.6}),
\begin{align}\label{W1}
W(z,s)=(T-t)^{1/(p-1)} y(x,t)
\end{align}
with
\begin{align}x-a=(T-t)^{1/2}z,\ T-t=e^{-s}, \end{align}
then $W(z,s)$ satisfies the following equation in $(\Omega-a) e^{s/2}\times[s_0(=-\log T), \infty)$,
\begin{align}\label{1.11}
\partial_s W=\partial_z^2 W-\frac{1}{2}z\partial_z W-\frac{1}{p-1}W
+|W|^{p-1}W+(\chi_\omega-1)|W|^{p-1}W.
\end{align}

For all $z\in \mathbb{R}^1$, we define
\begin{align}\label{w}w(z,s)=\left\{\begin{array}{ll} W(z,s)\chi(z,s), \ z\in (\Omega-a) e^{s/2},\\0, \ \ \ \ \ \ \ \ \ \ \ \ \ \ \ \ \ \ \mbox{otherwise},
\end{array}\right.\end{align}
with
\begin{align}\label{chi}\chi(z,s)=\chi_0\Big(\frac{ze^{-s/2}}{\varepsilon_0}\Big),\end{align}
where $\varepsilon_0>0$ with $(a-2\varepsilon_0,a+2\varepsilon_0)\subset \omega$ will be fixed small enough later.
Hence, $w(z,s)=0$, $|z|\geq 2\varepsilon_0 e^{s/2}$, from which and the fact $(a-2\varepsilon_0,a+2\varepsilon_0)\subset \omega$,  it holds
 that if the initial data $y_0\in C_0^\infty(\Omega)$ of support in $\omega$ and $y$ exists on $[0,T)$, then by the internal regularity of heat
   equation, $w$ is in the space
$C^\infty(\mathbb{R}\times[s_0,+\infty))$ with $s_0=-\log T$.

Then, we multiply equation (\ref{1.11}) by $\chi(z,s)$ and
get
\begin{align}\partial_s w=\partial_z^2 w-\frac{1}{2}z\partial_z w-\frac{1}{p-1}w
+|w|^{p-1}w+N(z,s),\ \forall\ z\in \mathbb{R},\ s_0\geq-\log T,\nonumber\end{align}
where
\begin{align}N(z,s)=\left\{\begin{array}{ll} W\partial_s\chi-2\partial_zW\partial_z\chi-W\partial_z^2\chi+\frac{1}{2}zW\partial_z\chi\\+|W|^{p-1}W(\chi-\chi^p),\  \ \ \ \ \ \ z\in (\Omega-a) e^{s/2},\ \ s_0\geq-\log T,\\0,\ \ \ \ \ \ \ \ \ \ \ \ \ \ \ \ \ \ \ \ \ \ \ \ \ \ \ \ \ \ \   z\in \mathbb{R}\setminus (\Omega-a)e^{s/2}, \ s_0\geq-\log T .
\end{array}\right.\end{align}
Here, we have use the fact that when $(a-2\varepsilon_0,a+2\varepsilon_0)\subset \omega$, $\chi_\omega\chi=\chi$.

Let \begin{align}\label{w1}w=\varphi+q,\end{align}
where  \begin{align}\label{w2}\varphi=f(\frac{z}{\sqrt{s}})+\frac{\kappa}{2ps}\end{align}
with $\kappa=(p-1)^{-\frac{1}{p-1}}$ and $f$ is defined in (\ref{f}). Then
$q$ is a solution to the following equation in $\mathbb{R}\times[s_0(=-\log T), \infty)$,
\begin{align}\label{2.20}
\partial_s q=(\mathcal{L}+V)q+B(z,s)+R(z,s)+N(z,s),
\end{align}
where
\begin{align}\label{}&
\mathcal{L}=\partial_z^2-\frac{1}{2}z\partial_z +1,\
V=p\varphi^{p-1}-\frac{p}{p-1},\\
\label{}&B(z,s)=|(\varphi+q)|^{p-1}(\varphi+q)-\varphi^p-p\varphi^{p-1}q,\\
\label{}&R(z,s)=\partial_z^2\varphi-\frac{1}{2}z\partial_z \varphi-\frac{1}{p-1}\varphi+\varphi^p-\partial_s \varphi,\\
\label{}&N(z,s)=H+\partial_z G(z,s),\\
\label{}&H(z,s)=W(\partial_s \chi+\partial_z^2 \chi+\frac{1}{2}z\partial_z\chi)+|W|^{p-1}W(\chi-\chi^p),\end{align}
and
\begin{align}\label{}&G(z,s)=-2\partial_z\chi W,\ \label{}\partial_zG(z,s)=-2\partial_z^2\chi W-2\partial_z\chi\partial_z W.\end{align}

We give a decomposition of the solution according to the spectrum of $\mathcal{L}$.
The operator $\mathcal{L}$ is self-adjoint on $\mathcal{D}(\mathcal{L})\subset L_\mu^2(\mathbb{R})$ with
\begin{align} \mu(z)=\frac{e^{-\frac{z^2}{4}}}{\sqrt{4\pi}},\end{align}
and $$L_\mu^2(\mathbb{R})=\Big\{v\in L_{loc}^2(\mathbb{R});\ \|v\|_{L_\mu^2}^2=\int |v|^2\mu(z)dz<+\infty\Big\}.$$
The spectrum of $\mathcal{L}$ is explicitly given by
$$\mbox{spec}(\mathcal{L})=\{1-\frac{m}{2};\ m\in \mathbb{N}\}.$$
 All the eigenvalues are simple. For $1-\frac{m}{2}$ corresponds the eigenfunction
\begin{align}
h_m(z)=\sum\limits_{n=0}^{[\frac{m}{2}]}\frac{m!}{n!(m-2n)!}(-1)^nz^{m-2n},
\end{align}
$h_m$ satisfies
$$\int h_nh_m\mu dz=2^n n!\delta_{nm}.$$
We will note also $k_m=h_m/\|h_m\|_{L_\mu^2}^2.$

First, let us introduce
\begin{align}\label{chi1}
\chi_1(z,s)=\chi_0\Big(\frac{|z|}{K_0\sqrt{s}}\Big),
\end{align}
where $\chi_0$ is defined in (\ref{1.7}), $K_0\geq 1$ will be chosen large enough.

We write $q=q_e+q_b$, where
\begin{align}
q_b=q\chi_1,\ q_e=q(1-\chi_1).
\end{align}
Then,
$$\mbox{supp} q_b(s)\subset B(0,2K_0\sqrt{s}),\ \mbox{supp} q_e(s)\subset \mathbb{R}\setminus B(0,K_0\sqrt{s}).$$
Second, we decompose $q_b$ as follows,
\begin{align}
q_b(z,s)=\sum\limits_{m=0}^2 q_m(s)h_m(z)+q_{-}(z,s),
\end{align}
where $q_m$ is the projection of $q_b$ on $h_m$, $q_{-}(z,s)=P_{-}(q_b)$, and $P_{-}$ is the projection on $\{h_i;\ i\geq 3\}$ the negative subspace of the operator $\mathcal{L}$.

In order to  reduce the infinite dimensional problem  proposed in Theorem \ref{Pros1} into finite
 dimensional one, we need the following definition as in \cite{Zaag2}.

\begin{definition}\label{s} For all $K_0>0$, $\varepsilon_0>0$, $A>0$, $0<\eta_0\leq 1$ and $T>0$, we define for all $t\in [0,T)$, the set
$S^*(K_0,\varepsilon_0,A,\eta_0,T,t)$ as being the set of all functions $y\in L^{\infty}(\Omega\times[0,T))$ satisfying

(i) Estimate in $\mathcal{R}_1$: $q(s)\in V_{K_0,A}(s)$,  where $s=-\log(T-t)$, $q(s)$ is defined in (\ref{W1}), (\ref{w}), (\ref{w1}), (\ref{w2})
and $V_{K_0,A}(s)$ is the set of all functions $r\in L^{\infty}(\mathbb{R}\times[-\log T,+\infty))$ such that
\begin{align}\left\{\begin{array}{ll} |r_m(s)|\leq As^{-2} (m=0, 1),\ |r_2(s)|\leq A^2s^{-2}\log s,\\|r_{-}(z,s)|\leq As^{-2}(1+|z|^3),\ \ \ \ |r_e(z,s)|\leq A^2s^{-1/2},
\end{array}\right.\end{align}
where
\begin{align}\left\{\begin{array}{ll} r_e(z,s)=(1-\chi_1(z,s))r(z,s),\ r_{-}(z,s)=P_{-}(\chi_1r),\\
\mbox{for}\ m\in \mathbb{N},\ r_m(s)=\int d\mu k_m(z)\chi_1(z,s)r(z,s).
\end{array}\right.\end{align}

(ii)  Estimate in $\mathcal{R}_2$: $\mbox{For all}\ \ x\in \Omega$ with $|x-a|\geq\frac{\varepsilon_0}{2}$, $|y(x,t)|\leq \eta_0$.
\end{definition}

For simplicity, we may write $S^*(t)$ instead of $S^*(K_0,\varepsilon_0,A,\eta_0,T,t)$.

\label{W}For all $K_0>0$, $\varepsilon_0>0$ and $A\geq 1$, there exists $\overline{s}(K_0,\varepsilon_0,A)>0$, such that if $s\geq  s_0\geq \overline{s}$, $0<\eta_0\leq 1$ and $y(t)\in S^*(t)$,
where $t=T-e^{-s}$, then we have
\begin{align}\label{W4}\left\{\begin{array}{ll}
\|q(s)\|_{L^\infty}\leq C A^2s^{-1/2},\ |q(z,s)|\leq CA^2\frac{\log s}{s^2}(1+|z|^3),\ \forall z\in \mathbb{R}, \\\
\|W(s)\|_{L^\infty}\leq \kappa+2.
\end{array}\right.\end{align}
Here and throughout the paper, $C$ is a constant number, which may be different in different context.

\subsection{Preparation of initial data}
Given $a\in \omega$, we consider initial data for system (\ref{1.6}) defined for all $x\in \Omega$ by
\begin{align}\label{initial data}
y_0(x,d_0,d_1)=T^{-\frac{1}{p-1}}\Big\{\varphi(z,s_0)\chi(8z,s_0)+\frac{A}{s_0^2}(d_0+d_1z)\chi_1(2z,s_0)\Big\},
\end{align}
where $T>0$ will be sufficiently small,
$s_0=-\log T$, $z=\frac{x-a}{\sqrt{T}}$, $\chi$ is defined in (\ref{chi}) and $\chi_1$ is defined in (\ref{chi1}).

Since $a\in\omega$, there exists $\delta_0\in (0,1)$ such that for any $\varepsilon_0$ with $\varepsilon_0\in(0,\delta_0]$, $(a-2\varepsilon_0,a+2\varepsilon_0)\subset \omega$. Then
  for any
$K_0>0$ and any $\varepsilon_0\in(0,\delta_0]$, there exits $\widetilde{s}_0(K_0,\varepsilon_0)$ such that when $s_0\geq  \widetilde{s}_0(K_0,\varepsilon_0)$, $K_0\sqrt{s_0}<\varepsilon_0e^{s_0/2}$
and $e^{s_0/2}/4<e^{s_0/2}$. Thus, $y_0(\cdot,d_0,d_1)\in C_0^\infty(\Omega)$ of support in $\omega$,
$\chi(8z,s_0)\chi(z,s_0)=\chi(8z,s_0)$ and
$\chi_1(2z,s_0)\chi(z,s_0)=\chi_1(2z,s_0)$ for $s_0\geq  \widetilde{s}_0(K_0,\varepsilon_0)$.

Then, by (\ref{W1}), (\ref{w}), (\ref{w1}) and (\ref{w2}), when $s_0\geq \widetilde{s}_0(K_0,\varepsilon_0)$,
\begin{align}&q(z,s_0)\nonumber\\
=&\left\{\begin{array}{ll}\varphi(z,s_0)(\chi(8z,s_0)-1)+\frac{A}{s_0^2}(d_0+d_1z)\chi_1(2z,s_0),\
z\in  (\Omega-a) e^{s_0/2},
\\q(z,s_0)=-\varphi(z,s_0),\ z\in \mathbb{R}\setminus (\Omega-a) e^{s_0/2}.\end{array}\right.\label{initial data3}\end{align}

In this subsection, we will prove that there exists $(d_0,d_1)\in \mathbb{R}^2$ such that the solution $y$ of (\ref{1.6}) with initial date $y_0$ defined in (\ref{initial data}) satisfies $y(t)\in S^*(t)$, $t\in[0,T)$, which will imply Theorem \ref{Pros1} (see section 2.4). First, we find a set $D_T\subset\mathbb{R}^2$ such that $y(0)\in S^*(0)$. More precisely, we have  the following lemma.

\begin{lemma} (Reduction for initial data) \label{initial data1} For any  $K_{0}>0$, for any  $\varepsilon_0\in(0,\delta_0]$ and for any $A\geq 1$, there exists  $s_{0,1}(K_0,\varepsilon_0,A)>0$ such that for
$s_0\geq s_{0,1}(K_0,\varepsilon_0,A)$ and $0<\eta_0\leq 1$:

If initial data for equation (\ref{1.6}) are given by  (\ref{initial data}), then there exists a rectangle
\begin{align}D_T\subset [-2,2]^2,\end{align}
such that

(i) for all $(d_0,d_1)\in D_T$, we have
$$y_0(\cdot,d_0,d_1)\in S^*(0);$$

(ii) for all $(d_0,d_1)\in D_T$ and $x\in \Omega$ with $|x-a|\geq \varepsilon_0/4$,
$$y_0(\cdot,d_0,d_1)=0.$$
\end{lemma}

Proof.
Let $K_0>0$,  $0<\varepsilon_0\leq \delta_0$ and $0<\eta_0\leq 1$.   We shall prove Lemma \ref{initial data1} in the following three steps.

Step 1. We will prove that for each $A\geq1$, there exists $s_0^1(K_0,\varepsilon_0,A)>0$ such that
for each $s_0\geq s_0^1(K_0,\varepsilon_0,A)$, if $(d_0,d_1)$ is such that $(q_0,q_1)(s_0)\in [-\frac{A}{s_0^2}, \frac{A}{s_0^2}]^2$, then
$$|d_0|\leq2,\ |d_1|\leq 2,$$
$$|q_2(s_0)|\leq (\log s_0)s_0^{-2},$$
$$|q_{-}(z,s_0)|\leq (1+|z|^3) s_0^{-2},$$
$$|q_e(\cdot,s_0)|\leq s_0^{-1/2}.$$

 Let $(d_0,d_1)\in \mathbb{R}^2$. When $z\in (\Omega-a) e^{s_0/2}$, then  by (\ref{initial data3}), if $s_0\geq  \widetilde{s}_0(K_0,\varepsilon_0)$,
we can write initial data as
$$q(z,s_0)=q^0(z,s_0)+q^1(z,s_0)+q^2(z,s_0)+q^3(z,s_0),$$ where
 $q^0(z,s_0)=\Big(p-1+\frac{(p-1)^2}{4p}\frac{z^2}{s_0}\Big)^{-\frac{1}{p-1}}(\chi(8z,s_0)-1)$,
 $q^1(z,s_0)=\frac{\kappa}{2ps_0}(\chi(8z,s_0)-1)$, $q^2(z,s_0)=\frac{A}{s_0^2}d_0
 \chi_1(2z,s_0)$ and  $q^3(z,s_0)=\frac{A}{s_0^2}
 d_1z\chi_1(2z,s_0)$; If $z\in \mathbb{R}\setminus (\Omega-a) e^{s_0/2}$,
$q=-\varphi=-\Big(p-1+\frac{(p-1)^2}{4p}\frac{z^2}{s_0}\Big)^{-\frac{1}{p-1}}-\frac{\kappa}{2ps_0}$.

a) Estimate for $d_1$, $d_2$ and $q_2(s_0)$.

By the definition of $\chi$ in (\ref{chi}) and the definition of $\chi_1$  in (\ref{chi1}), there exists $s_0^2(K_0,\varepsilon_0)$ such that
if $s_0\geq s_0^2(K_0,\varepsilon_0)$, then  $2K_0\sqrt{s_0}<\varepsilon_0 e^{s_0/2}/8$.
Thus, for $s_0\geq s_0^2(K_0,\varepsilon_0)$, \begin{align}\chi_1(z,s_0)(\chi(8z,s_0)-1)=0.\end{align}
Hence, if $s_0\geq s_0^2(K_0,\varepsilon_0)$, $q^0(z,s_0)\chi_1(z,s_0)=0$, $q^1(z,s_0)\chi_1(z,s_0)=0$, $\forall z\in (\Omega-a) e^{s_0/2},
$ and $$\int_{(\Omega-a) e^{s_0/2}}
 q^0k_m\chi_1(z,s_0)d\mu=0,\ \int_{(\Omega-a) e^{s_0/2}}
 q^1k_m\chi_1(z,s_0)d\mu=0, \ m=0,1,2.$$

Since $a\in \omega$, it holds that $\{z;\ |z|\leq K_0\sqrt{s_0}\}\subset(\Omega-a) e^{s_0/2}$ if $s_0\geq s_0^3(K_0,\varepsilon_0)$ with $s_0^3$ large enough. Thus, \begin{align}\int_{(\Omega-a) e^{s_0/2}}
 q^2k_0\chi_1(z,s_0)d\mu= d_0\frac{A}{s_0^2}\int_{(\Omega-a) e^{s_0/2}}
 \chi_1(2z,s_0)\chi_1(z,s_0)d\mu\nonumber\\ =d_0\frac{A}{s_0^2}\int_{(\Omega-a) e^{s_0/2}}
 \chi_1(2z,s_0)d\mu=d_0\frac{A}{s_0^2}\int
 \chi_1(2z,s_0)d\mu=:d_0\frac{A}{s_0^2} b_0(s_0),\nonumber\end{align}
 where $|b_0(s_0)|\leq \int \frac{d\mu}{\sqrt{4\pi}}= 1$.
 Moreover,  if $s_0\geq s_0^3(K_0,\varepsilon_0)$,  \begin{align}\Big|\int_{(\Omega-a) e^{s_0/2}}
 q^2k_0\chi_1(z,s_0)d\mu\Big|=\Big|d_0\frac{A}{s_0^2}\int
 (1+(\chi_1(2z,s_0)-1))d\mu\Big|\leq  |d_0|\frac{A}{s_0^2}(1+Ce^{-Cs_0}).\nonumber\end{align}
In the above estimate, we have used that when
$$|z|\geq  \frac{1}{2}K_0\sqrt{s_0},\ \exp{\Big(-\frac{|z|^2}{4}\Big)}\leq \exp{\Big(-\frac{|z|^2}{8}\Big)}\exp{\Big(-\frac{C(K_0)s_0}{8}\Big)}.$$
When $s_0\geq s_0^3(K_0,\varepsilon_0)$, we also have the following equalities
$$\int_{(\Omega-a) e^{s_0/2}}
 q^2k_1\chi_1(z,s_0)d\mu= d_0\frac{A}{s_0^2}\int_{-K_0\sqrt{s_0}}^{K_0\sqrt{s_0}}
 \chi_1(2z,s_0)\chi_1(z,s_0)\frac{z}{2}d\mu=0,$$
 and
\begin{align}\label{equality}&\Big|\int_{(\Omega-a) e^{s_0/2}}
 q^2k_2\chi_1(z,s_0)d\mu|
 \nonumber\\= &\Big|d_0\frac{A}{s_0^2}\int
 (\chi_1(2z,s_0)-1)\frac{z^2-2}{8}d\mu\Big|
 \leq C |d_0|\frac{A}{s_0^2}e^{-Cs_0}.\end{align}
Similarly,  when $s_0\geq s_0^3(K_0,\varepsilon_0)$,
$$\int_{(\Omega-a) e^{s_0/2}}
 q^3k_0\chi_1(z,s_0)d\mu=d_1\frac{A}{s_0^2}\int_{(\Omega-a) e^{s_0/2}}
z\chi_1(2z,s_0)\chi_1(z,s_0)d\mu=0,$$
\begin{align}&\int_{(\Omega-a) e^{s_0/2}}
 q^3k_1\chi_1(z,s_0)d\mu= d_1\frac{A}{s_0^2}\int_{(\Omega-a) e^{s_0/2}} \frac{z^2}{2}(1+
 \chi_1(2z,s_0)-1)d\mu=: d_1\frac{A}{s_0^2} b_1(s_0),\nonumber\end{align}
where $|b_1(s_0)|\leq\int\frac{z^2}{2}d\mu= 1$,
  \begin{align}\Big|&\int_{(\Omega-a) e^{s_0/2}}
 q^3k_1\chi_1(z,s_0)d\mu\Big|=\Big|d_1\frac{A}{s_0^2}\Big(1+\int
 (\chi_1(2z,s_0)-1)d\mu\Big)\Big|\nonumber\\\leq & |d_1|\frac{A}{s_0^2}(1+Ce^{-Cs_0}),\nonumber\end{align}
 and
$$\int_{(\Omega-a) e^{s_0/2}}
 q^3k_2\chi_1(z,s_0)d\mu=d_1\frac{A}{s_0^2}\int_{-K_0\sqrt{s_0}}^{K_0\sqrt{s_0}}
 z\chi_1(2z,s_0)\chi_1(z,s_0)\frac{z^2-2}{8}d\mu=0.$$

On the other hand, since $a\in \omega$, when $s_0\geq s_0^4(K_0,\varepsilon_0)$ with $s_0^4$ large enough, and $z\in \mathbb{R}\setminus (\Omega-a) e^{s_0/2}$, we have
$$|z|\geq 2K_0\sqrt{s_0},$$ which implies $q(z,s_0)\chi_1(z,s_0)=0$, $z\in \mathbb{R}\setminus (\Omega-a) e^{s_0/2}$, and
$$\int_{\mathbb{R}\setminus (\Omega-a) e^{s_0/2}}
 q(z,s_0)k_m\chi_1(z,s_0)d\mu=0,\ m=0,1,2.$$
 Hence, it holds that when $s_0\geq s_0^5(K_0,\varepsilon_0)$ with $s_0^5$ large enough,  $q_0(s_0)=d_0\frac{A}{s_0^2} b_0(s_0)$, $q_1(s_0)=d_1\frac{A}{s_0^2} b_1(s_0)$,
 where $|{b}_0(s_0)|\leq 1$, $|{b}_1(s_0)|\leq 1$.
Furthermore, since $${\int 1d\mu}= 1,\ {\int_{-K_0\sqrt{s_0}/2}^{K_0\sqrt{s_0}/2}1d\mu}\geq 1/2$$
 for $s_0\geq s_0^6(K_0,\varepsilon_0)$ with $s_0^6$ large enough. If $(d_0,d_1)$ is chosen such that $(q_0,q_1)(s_0)\in [-\frac{A}{s_0^2}, \frac{A}{s_0^2}]^2$, then we have that if $s_0\geq s_0^7(K_0,\varepsilon_0)$ with $s_0^7$ large enough, then
 $$|d_0|\leq \frac{1}{ b_0(s_0)}=\frac{1}{\int_{(\Omega-a) e^{s_0/2} }
 \chi_1(2z,s_0)d\mu}\leq \frac{1}{\int_{-K_0\sqrt{s_0}/2}^{K_0\sqrt{s_0}/2}1d\mu}\leq 2.
 $$

Similarly, since $${\int \frac{z^2}{2}d\mu}= 1,$$  we have that when $s_0\geq s_0^8(K_0,\varepsilon_0)$ with $s_0^8$ large enough, if $(d_0,d_1)$ is chosen such that $(q_0,q_1)(s_0)\in [-\frac{A}{s_0^2}, \frac{A}{s_0^2}]^2$, then
  $$|d_1|\leq \frac{1}{ b_1(s_0)}=\frac{1}{\int_{(\Omega-a) e^{s_0/2} }
 \chi_1(2z,s_0)\frac{z^2}{2}d\mu}\leq \frac{1}{\int_{-K_0\sqrt{s_0}/2}^{K_0\sqrt{s_0}/2}\frac{z^2}{2}d\mu}\leq 2.
 $$

On the other hand, by (\ref{equality}), if $s_0\geq s_0^9(K_0,\varepsilon_0)$ with $s_0^9$ large enough, then
$$|q_2(s_0)|\leq \frac{\log s_0}{s_0^2}.$$

b) Estimate for $q_{-}(z,s_0)$.

Since $$|q_{-}(z,s_0)|=\Big|q(z,s_0)\chi_1(z,s_0)-\sum\limits_{m=0}^2 q_m(s)h_m(z)\Big|,\ z\in \mathbb{R},$$
it holds that when $s_0 \geq s_0^{10}(K_0,\varepsilon_0,A)$ with $s_0^{10}$ large enough,
\begin{align}|q_{-}(z,s_0)|\leq \nonumber& C\frac{A}{s_0^2}\Big(\frac{|z|^3}{{s_0^{3/2}}}+e^{-Cs_0}(1+|z|^2)\Big)+C\frac{A}{s_0^2}\Big(\frac{|z|^3}{{s_0}}+e^{-Cs_0}|z|\Big)\nonumber\\
\leq &Cs_0^{-2}(1+|z|^3),\ z\in (\Omega-a) e^{s_0/2},\nonumber\end{align}
and
\begin{align}&|q_{-}(z,s_0)|\leq C|d_0|\frac{A}{s_0^2}+C|d_1|\frac{A}{s_0^2}|z|+ C(|d_0|+|d_1|)\frac{A}{s_0^2}e^{-Cs_0}(1+|z|^2)\nonumber\\\leq
&C{s_0^{-2}}(1+|z|^3),\ z\in \mathbb{R}\setminus (\Omega-a) e^{s_0/2}.\end{align}
Here we have used that when $z\in \mathbb{R}\setminus (\Omega-a) e^{s_0/2}$ with  $s_0$ large enough, $|z|\geq{s_0}$.

c) Estimate for $q_{e}(z,s_0)$.

Since \begin{align}\chi_1(2z,s_0)(1-\chi_1(z,s_0))=0,\ z\in \mathbb{R},\end{align}
it follows that $$(1-\chi_1)q^2(z,s_0)=0,\ (1-\chi_1)q^3(z,s_0)=0,\ z\in  (\Omega-a) e^{s_0/2}.$$

On the other hand, if
$s_0 \geq s_0^{11}$ with $s_0^{11}$ large enough,
$$|(1-\chi_1)q^1(z,s_0)|\leq \frac{C}{s_0}\leq s_0^{-1/2}/3,\ z\in  (\Omega-a) e^{s_0/2}.$$

When $s_0\geq s_0^{12}(K_0,\varepsilon_0,A)$ with $s_0^{12}$ large enough,
\begin{align}|(1-\chi_1)q^0(z,s_0)|
=\Big|\Big(p-1+\frac{(p-1)^2}{4p}\frac{z^2}{s_0}\Big)^{-\frac{1}{p-1}}(\chi(8z,s_0)-1)\Big|,\ z\in  (\Omega-a) e^{s_0/2}.\nonumber\end{align}
This implies that when $s_0\geq s_0^{12}$ and $\frac{|z|}{\varepsilon_0 e^{s_0/2}/8}\geq 1$, $|(1-\chi_1)q^0| \neq 0$. Thus, When $s_0\geq s_0^{13}(K_0,\varepsilon_0,A)$ with $s_0^{13}$ large enough,
\begin{align}|(1-\chi_1)q^0(z,s_0)|\leq &C(\frac{z^2}{s_0})^{-\frac{1}{p-1}}\frac{|z|^{\frac{2}{p-1}}}{(\varepsilon_0 e^{s_0/2}/8)^{\frac{2}{p-1}}}\nonumber\\\leq& C(\frac{1}{s_0})^{-\frac{1}{p-1}}(\frac{1}{\varepsilon_0} e^{-s_0/2})^{\frac{2}{p-1}}\leq \frac{C}{s_0}\leq s_0^{-1/2}/3,\ z\in  (\Omega-a) e^{s_0/2}.\nonumber\end{align}
When $z\in \mathbb{R}\setminus (\Omega-a) e^{s_0/2}$,
\begin{align}|(1-\chi_1)q|\leq (p-1+\frac{ (C e^{s_0/2})^2}{s_0})^{-\frac{1}{p-1}}+\frac{\kappa}{2ps_0}\leq s_0^{-1/2}/3.\nonumber\end{align}
if $s_0\geq s_0^{14}(K_0,\varepsilon_0,A)$ with $s_0^{14}$ large enough,.

Hence, we have $$|q_e(z,s_0)|\leq s_0^{-1/2},\ z\in \mathbb{R}.$$

Step 2. We will prove that  for any $A\geq 1$ and
for any $s_0\geq s_0^{15}(K_0,\varepsilon_0,A)$,  there exists a set $\mathcal{D}_{s_0}\subset \mathbb{R}^2$
topologically equivalent to a square with the following property:
$$q(d_0,d_1,s_0)\in V_{K_0,A}(s_0)\ \mbox{if and only if}\ (d_0,d_1)\in \mathcal{D}_{s_0}.$$

By step 1,  for each $A\geq 1$, there exists $s_0^1(K_0,\varepsilon_0,A)>0$ such that
for each $s_0\geq s_0^1(K_0,\varepsilon_0,A)$, we have the following equivalence,
$$q(s_0)\in  V_{K_0,A}(s_0)\ \mbox{if and only if } (q_0,q_1)(s_0)\in [-\frac{A}{s_0^2}, \frac{A}{s_0^2}]^2.$$
Then it is enough to prove that there exists a set  $\mathcal{D}_{s_0}\subset \mathbb{R}^2$ topologically equivalent to a square satisfying
$$(q_0,q_1)(s_0)\in [-\frac{A}{s_0^2}, \frac{A}{s_0^2}]^2\ \mbox{if and only if}\ (d_0,d_1)\in \mathcal{D}_{s_0}.$$
We are sufficient to take
$$\mathcal{D}_{s_0}=\Big[-\frac{1}{b_0(s_0)},\frac{1}{b_0(s_0)}\Big]\times\Big[-\frac{1}{b_1(s_0)},\frac{1}{b_1(s_0)}\Big]\subset [-2,2]\times[-2,2],$$
 when $s_0\geq s_0^{15}(K_0,\varepsilon_0,A)$.

 Step 3. Suppose that $x\in \Omega$ with $|x-a|\geq \varepsilon_0/4$, we have $|z|\geq e^{s_0/2}\varepsilon_0/4$.
In this case, there exists $s_0^{16}(K_0,\varepsilon_0)>0$ such that when $s_0\geq s_0^{16}(K_0,\varepsilon_0)$, $\varepsilon_0 e^{s_0/2}/4\geq K_0\sqrt{s_0}$, $\chi(8z,s_0)=0$ and $\chi_1(2z,s_0)=0$. Thus, by the definition of $y_0(\cdot,d_0,d_1)$ in (\ref{initial data}), we have that when $s_0\geq s_0^{16}(K_0,\varepsilon_0)$, for any $(d_0,d_1)\in \mathcal{D}_{T}$ and  for any $x$ with $|x-a|\geq \varepsilon_0/4$,
$$y_0(x,d_0,d_1)=0.$$

By step 1, step 2, and step 3, we can get Lemma \ref{initial data1}.\ \ \ \#

\subsection{Reduction to a finite-dimensional problem}

Let us consider $(d_0,d_1)\in D_T$, and $s_0=-\log T\geq s_{0,1}$, defined in Lemma \ref{initial data1}. By classical theory of parabolic equations, we can define a
maximal solution $y$ to equation (\ref{1.6}) with initial data (\ref{initial data}), and a maximal time $t_*(d_0,d_1)\in [0,T]$ such that
\begin{align}
\forall t\in[0,t_*),\ y(t)\in S^*(t).
\end{align}

(1) either $t_*=T$,

(2) or $t_*<T$ and from continuity, $y(t_*)\in \partial S^*(t_*)$,
in the sense that when $t=t_*$, one `$\leq$' symbol in the definition of $S^*(t_*)$ is replaced by the symbol ``$=$".

Our aim is to show that for all $A$ large and $T$ small enough, one can find $(d_0,d_1)\in D_T$, $t_*(d_0,d_1)=T$,

We argue by contradiction to prove Theorem \ref{Pros1}, and assume that for all $(d_0,d_1)\in D_T$, $t_*(d_0,d_1)<T$. We will reduce
the problem of controlling all the components of $y$ in $S^*(t)$ to a problem of controlling $(q_0,q_1)(s)$.

\subsubsection{Priori estimates}

By the definition of $S^*(t)$, there are two different types of estimates in the regions $\mathcal{R}_1$ and $\mathcal{R}_2$. Thus, in this subsection, we give
two different priori estimates in $\mathcal{R}_1$ and $\mathcal{R}_2$, respectively.

{\bf Part 1: Estimates in $\mathcal{R}_1$}

The following proposition gives a priori estimate in $\mathcal{R}_1$.
\begin{proposition}\label{qb}There exists $K_{0}^1>0$ such that for any $K_0\geq K_0^1$ and for any $\varepsilon_0\in(0,\delta_0]$,  there exists $A_1({K_0,\varepsilon_0})$
such that for any $A\geq A_1$ and $0<\eta_0\leq 1$, there exists $\widehat{s}_1(K_0,\varepsilon_0,A)$ such that for all  $s_0\geq \widehat{s}_1(K_0,\varepsilon_0,A)$, for any solution $q$ of (\ref{2.20}),
we have the following property: if $(d_0,d_1)$ is chosen so that $(q_0(s_0),q_1(s_0))\in [-\frac{A}{s_0^2},\frac{A}{s_0^2}]^2$, and
            if for $s_1\geq s_0$, we have $\forall s\in[s_0,s_1]$, $y(t)\in S^*(t)$ with $t=T-e^{-s}$, then $\forall s\in[s_0,s_1]$,
            \begin{align}\label{p2}
               |q_2(s)|\leq {A^2}s^{-2}\log s-s^{-3},\nonumber\\
             |q_{-}(z,s)|\leq {\frac{A}{2}}(1+|z|^3)s^{-2},\nonumber\\
            \|q_e(s)\|_{L^\infty}\leq \frac{A^2}{2\sqrt{s}}. \end{align}
Here $\delta_0$ is defined in section 2.2.
\end{proposition}

Let us first write equation (\ref{2.20}) in its Duhamel formulation,

\begin{align}\label{q}
q(s)=&K(s,\sigma)q(\sigma)+\int_\sigma^sd\tau K(s,\tau)B(q(\tau))+\int_\sigma^sd\tau K(s,\tau)R(\tau)\nonumber\\&+\int_\sigma^sd\tau K(s,\tau)(H+\partial_zG)(\tau),
\end{align}
where $K$ is the fundamental solution of the operator $\mathcal{L}+V$ defined for each $s_0>0$, for each $s_1\geq s_0$ by
\begin{align}\label{1.1}
&\partial_{s_1} K(s_1,s_0)=(\mathcal{L}+V)K(s_1,s_0),\nonumber\\
&K(s_0,s_0)=I.
\end{align}
 We write $q=\alpha+\beta+\gamma+\delta+\overline{\delta}$ with
\begin{align}\label{1.1}
\alpha(s)=K(s,\sigma)q(\sigma),\ \beta(s)=\int_\sigma^s d\tau K(s,\tau)B(q(\tau)), \nonumber\\
\gamma(s)=\int_\sigma^s d\tau K(s,\tau)R(\tau), \ \delta(s)=\int_\sigma^s d\tau K(s,\tau)H(\tau), \\
\overline{\delta}(s)=\int_\sigma^s d\tau K(s,\tau)\partial_zG(\tau),\nonumber
\end{align}
where for a function $F(z,\tau)$,  $K(s,\tau)F(\tau)$ is defined by
$$K(s,\tau)F(\tau)=\int dx K(s,\tau,z,x)F(x,\tau).$$

In order to prove Proposition \ref{qb}, we need the following lemma.

            \begin{lemma}\label{BK} a) $\forall s\geq \tau\geq 1$ with $s\leq 2\tau$, $\forall z,x\in \mathbb{R}$,
            $|K(s,\tau,z,x)|\leq Ce^{(s-\tau)\mathcal{\mathcal{L}}}(z,x)$ with $$e^{\theta\mathcal{L}}(z,x)=\frac{e^\theta}{\sqrt{4\pi(1-e^{-\theta})}}\exp\Big[-\frac{(ze^{-\theta/2}-x)^2}{4(1-e^{-\theta})}\Big].$$

           b) For all $s\geq \tau \geq 1$ with $s\leq 2\tau$, we have the following : For all
differentiable function $g\in L^\infty$, such that $zg\in L^\infty$,
\begin{align}&\|K(s,\tau)\partial_zg\|_{L^\infty}\nonumber\\\leq &Ce^{s-\tau}\Big\{\frac{\|g\|_{L^\infty}}{\sqrt{1-e^{-(s-\tau)}}}+\frac{s-\tau}{s}(1+s-\tau)\Big((1+e^{\frac{s-\tau}{2}})\|zg\|_{L^\infty}
+e^{\frac{(s-\tau)}{2}}\|g\|_{L^\infty}\Big)\Big\}.\nonumber\end{align}

            c) There exists $K_{01}>0$ such that for any $K_0\geq K_{01}$,
            for any $A'>0$, $A''>0$, $A'''>0$ and $\rho^*>0$, there exists $s_1(A', A'',A''',K_0,\rho^*)$ with the following property:
            $\forall s_0\geq s_1$, assume that for $\sigma\geq s_0$,
            \begin{align}\label{qsigma}
              |q_m(\sigma)|\leq A'\sigma^{-2},\ m=0,1,\\
               |q_2(\sigma)|\leq A''(\log \sigma)\sigma^{-2},\nonumber\\
             |q_{-}(z,\sigma)|\leq A'''(1+|z|^3)\sigma^{-2},\nonumber\\
            \|q_e(\sigma)\|_{L^\infty}\leq A''\sigma^{-1/2}\nonumber, \end{align}
           then, $\forall s\in[\sigma,\sigma+\rho^*]$,
           \begin{align}|\alpha_2(s)|\leq  A''(\log \sigma)s^{-2} + C\max\{A',A'''\}(s-\sigma)e^{(s-\sigma)}s^{-3}, \nonumber\end{align}
           $$|\alpha_-(z,s)|\leq C (A^{'''}s^{-2}e^{-(s-\sigma)/2}+A''e^{-(s-\sigma)^2} s^{-2})(1+|z|^3),$$
             $$\|\alpha_e(s)\|_{L^\infty}\leq Ce^{(s-\sigma)}A'''s^{-1/2}+C A''s^{-1/2}e^{-(s-\sigma)/p},$$
             where $$\alpha(s)=K(s,\sigma)q(\sigma)=\sum_{m=0}^{2}\alpha_m(s) h_m(z)+\alpha_{-}(z,s)+\alpha_e(z,s).$$

             d) $\forall \rho^*>0$, there exists $s_{2}(\rho^*)$ such that $\forall \sigma\geq s_{2}(\rho^*)$, $\forall s\in [\sigma,\sigma+\rho^*]$,
              \begin{align}|\gamma_2(s)|\leq  C(s-\sigma)s^{-3}, \nonumber\end{align}
           $$|\gamma_-(z,s)|\leq C (s-\sigma) s^{-2}(1+|z|^3),$$
            where
            $$\int_\sigma^s K(s,\tau)R(\tau)=\gamma(z,s)=\sum_{m=0}^{2}\gamma_m(s) h_m(z)+\gamma_{-}(z,s)+\gamma_e(z,s).$$
            \end{lemma}

            From a) of the above lemma, the following corollary holds (see Corrollary 3.14  in \cite{Zaag1}).
\begin{corollary} \label{corollary 3.1} $\forall s\geq \tau\geq 1$ with $s\leq 2\tau$,
\begin{align}
\Big|\int K(s,\tau,z,x)(1+|x|^m)dx\Big|\leq &C\int e^{(s-\tau)\mathcal{L}}(z,x)(1+|x|^m)dx\nonumber\\\leq &Ce^{(s-\tau)}(1+|z|^m), \ m\geq 1.\nonumber
\end{align}
\end{corollary}

Lemma \ref{BK} is similar to Lemma 3.13 in in \cite{Zaag1}.
 We shall take the  techniques in the proof of Lemma 3.13 in in \cite{Zaag1} to prove this lemma and just give a detailed proof for those different with  Lemma 3.13 in in \cite{Zaag1} (see Appendix).


We are now going to prove the following proposition which implies Proposition \ref{qb}.

\begin{proposition}\label{qa}There exists $K_{0}^1>0$ such that for any $K_0\geq K_0^1$,   any $\varepsilon_0\in(0,\delta_0]$ and any
$\widetilde{A}>0$,   there exists $\widetilde{A}_2(\widetilde{A},K_0,\varepsilon_0)>0$ such that for all $A\geq\widetilde{A}_2(\widetilde{A},K_0,\varepsilon_0)$,
there exists $\widehat{s}_2(\widetilde{A},A,K_0,\varepsilon_0)>0$ such that for any $s_0\geq \widehat{s}_2(\widetilde{A},A,K_0,\varepsilon_0)$,
for any solution of (\ref{2.20}),
we have the following property: if
   \begin{align}
              |q_m(s_0)|\leq As_0^{-2},\ m=0,1,\nonumber\\
               |q_2(s_0)|\leq \widetilde{A}s_0^{-2}\log s_0,\label{property}\\
             |q_{-}(z,s_0)|\leq \widetilde{A}(1+|z|^3)s_0^{-2},\nonumber\\
            \|q_e(s_0)\|_{L^\infty}\leq \widetilde{A}s_0^{-1/2}\nonumber, \end{align}
            if for $s_1\geq s_0$, we have $\forall s\in[s_0,s_1]$, $y(t)\in S^*(t)$ with $t=T-e^{-s}$, then $\forall s\in[s_0,s_1]$,
            \begin{align}\label{}
               |q_2(s)|\leq {A^2}s^{-2}\log s-s^{-3},\nonumber\\
             |q_{-}(z,s)|\leq {\frac{A}{2}}(1+|z|^3)s^{-2},\label{property1}\\
            \|q_e(s)\|_{L^\infty}\leq \frac{A^2}{2\sqrt{s}}\nonumber. \end{align}
\end{proposition}

{\bf Proposition \ref{qa} implies Proposition \ref{qb}}

Indeed, referring to step 1 Lemma \ref{initial data1}, we apply proposition \ref{qa}, with $\widetilde{A}=1$. This gives $\widetilde{A}_2(1,K_0,\varepsilon_0)>0$
and $\widehat{s}_2(1,A,K_0,\varepsilon_0)$. If we take $\widehat{s}_1(K_0,\varepsilon_0,A)=\max(\widehat{s}_2(1,A,K_0,\varepsilon_0),s_{0,1}(K_0,\varepsilon_0,A))$,
we can get Proposition \ref{qb}.\ \ \ \#

Proposition \ref{qa} can be derived from the following lemma.

\begin{lemma}\label{main}
There exists $K_{02}>0$ and $A_2>0$ such that for any $K_0\geq K_{02}$, $\varepsilon_0\in(0,\delta_0]$, $A\geq A_2$, $\widetilde{A}>0$ and $\rho^*>0$, there exists $s_2(A,\widetilde{A},\rho^*)>0$
with the following property: $\forall s_0\geq s_2(A,\widetilde{A},K_0,\varepsilon_0,\rho^*)$, $\forall \rho\leq \rho^*$, assume $\forall s\in[\sigma,\sigma+\rho]$  with $\sigma\geq s_0$, and $y(t)\in S^*(t)$ with $t=T-e^{-s}$.

I) Case $\sigma\geq s_0$: we have $\forall s\in[\sigma,\sigma+\rho]$,

i) (linear term)
\begin{align}|\alpha_2(s)|\leq  A^2(\log \sigma)s^{-2} + CA(s-\sigma)e^{(s-\sigma)}s^{-3}, \nonumber\end{align}
$$|\alpha_-(z,s)|\leq C (As^{-2}e^{-(s-\sigma)/2}+A^2e^{-(s-\sigma)^2} s^{-2})(1+|z|^3),$$
             $$\|\alpha_e(s)\|_{L^\infty}\leq Ce^{(s-\sigma)}As^{-1/2}+C A^2s^{-1/2}e^{-(s-\sigma)/p},$$
             where $$K(s,\sigma)q(\sigma)=\alpha(z,s)=\sum_{m=0}^{2}\alpha_m(s) h_m(z)+\alpha_{-}(z,s)+\alpha_e(z,s).$$
\end{lemma}

ii) nonlinear term
$$|\beta_2(s)|\leq\frac{s-\sigma}{s^{3+1/2}}, $$
$$|\beta_-(z,s)|\leq(s-\sigma)(1+|z|^3){s^{-2-\varepsilon}}, $$
$$\|\beta_e(s)\|_{L^\infty}\leq (s-\sigma)s^{-\frac{1}{2}-\varepsilon},$$
where $\varepsilon=\varepsilon_1(p)>0$, and
$$\int_\sigma^s d\tau  K(s,\tau)B(q(\tau))=\beta(z,s)=\sum_{m=0}^{2}\beta_m(s) h_m(z)+\beta{-}(z,s)+\beta_e(z,s).$$

iii) corrective term
  \begin{align}|\gamma_2(s)|\leq  (s-\sigma)Cs^{-3}, \nonumber\end{align}
           $$|\gamma_-(z,s)|\leq C (s-\sigma) s^{-2}(1+|z|^3),$$
           $$\|\gamma_e(s)\|_{L^\infty}\leq (s-\sigma)s^{-3/4},$$ where
$$\int_\sigma^s d\tau  K(s,\tau)R(\cdot,\tau)=\gamma(z,s)=\sum_{m=0}^{2}\gamma_m(s) h_m(z)+\gamma{-}(z,s)+\gamma_e(z,s).$$

iv) \begin{align}|\delta_m(s)|\leq C\frac{s-\sigma}{s^{3+1/2}},\nonumber\end{align}
\begin{align}|\delta_-(s)|\leq C\frac{s-\sigma}{s^{2+\varepsilon}}(1+|z|^3),\nonumber\end{align}
\begin{align}|\delta_e(s)|\leq C\frac{s-\sigma}{s^{1/2+\varepsilon}},\nonumber\end{align}
where
$$\int_\sigma^s d\tau  K(s,\tau)H(\cdot,\tau)=\delta(z,s)=\sum_{m=0}^{2}\delta_m(s) h_m(z)+\delta{-}(z,s)+\delta_e(z,s).$$

v) \begin{align}|\widetilde{\delta}_m(s)|\leq C\frac{(s-\sigma)+\sqrt{s-\sigma}}{s^{3+1/2}},\nonumber\end{align}
\begin{align}|\widetilde{\delta}_-(s)|\leq C\frac{(s-\sigma)+\sqrt{s-\sigma}}{s^{2+\varepsilon}}(1+|z|^3),\nonumber\end{align}
\begin{align}|\widetilde{\delta}_e(s)|\leq C\frac{(s-\sigma)+\sqrt{s-\sigma}}{s^{1/2+\varepsilon}},\nonumber\end{align}
where
$$\int_\sigma^s d\tau  K(s,\tau)\partial_zG(\cdot,\tau)=\widetilde{\delta}(z,s)=\sum_{m=0}^{2}\widetilde{\delta}_m(s) h_m(z)
+\widetilde{\delta}{-}(z,s)+\widetilde{\delta}_e(z,s).$$
II) Case $\sigma=s_0$. Assume in addition  $q(s_0)$ satisfies (\ref{property}). Then  $\forall s\in[s_0,s_0+\rho]$,

i) linear term
\begin{align}|\alpha_2(s)|\leq  \widetilde{A}(\log s_0)s^{-2} + C\max\{A,\widetilde{A}\}e^{(s-s_0)}(s-s_0)s^{-3}, \nonumber\end{align}
$$|\alpha_-(z,s)|\leq C \widetilde{A}s^{-2}(1+|z|^3),$$
             $$\|\alpha_e(s)\|_{L^\infty}\leq C(1+e^{(s-s_0)})\widetilde{A}s^{-1/2}.$$

Proof. We choose $s_0\geq \rho^*$ in all cases so that if $s_0\leq \sigma\leq \tau\leq s\leq \sigma+\rho^*$ and $\rho\leq \rho^*$,
we have $\sigma^{-1}\leq 2s^{-1}$ and $\tau^{-1}\leq 2s^{-1}$.

By Lemma \ref{BK}, we can use the similar techniques  in the proof of Lemma 3.12 in \cite{Zaag1} to prove Ii), IIi), Iii) and Iiii).

Iiv) Estimate of $\delta$.

Consider $s\in [\sigma,\sigma+\rho^*)$, we recall that $0<\eta_0\leq 1$. Since $y(t)\in S^*(t)$, where $t=T-e^{-s}$,
we see that if $s_0\geq s_2^{1}$ with $s_2^{1}$ large enough,
\begin{align}\label{}\|H(s)\|_{L^\infty}\leq \frac{C}{\varepsilon_0^2}e^{-\frac{s}{p-1}}.\nonumber\end{align}
In particular, if $\sigma\leq \tau\leq s\leq \sigma+\rho^*$,  then \begin{equation}\nonumber
\frac{1}{2\sigma}\leq \frac{1}{s}\leq \frac{1}{\tau}\le \frac 1 \sigma,\ \tau\geq \sigma\geq \frac{s}{2},
\end{equation}
\begin{align}\nonumber\|H(\tau)\|_{L^\infty}\leq \frac{C}{\varepsilon_0^2}e^{-\frac{\tau}{p-1}}\leq \frac{C}{\varepsilon_0^2}e^{-\frac{s}{2(p-1)}}.\end{align}
Hence, if $s_0\geq s_2^{2}(\rho^*)$ with $s_2^{2}$ large enough,
\begin{align}|\delta(z,s)|\leq&\int_\sigma^sd\tau\int |K(s,\tau,z,x)H(x,\tau)|dx\nonumber\\
\leq&\int_\sigma^sd\tau\int e^{(s-\tau)\mathcal{L}}(z,x)\frac{C}{\varepsilon_0^2}e^{-\frac{s}{2(p-1)}}dx\nonumber\\
\leq &\frac{C}{\varepsilon_0^2}e^{-\frac{s}{2(p-1)}} \int_\sigma^sd\tau e^{s-\tau}\leq \frac{C(\rho^*)}{\varepsilon_0^2}e^{-\frac{s}{2(p-1)}} (s-\sigma)\nonumber\\\leq &\frac{s-\sigma}{s^{3+1/2}}.\nonumber\end{align}
Thus, if  if $s_0\geq s_2^{3}(\varepsilon_0,\rho^*)$ with $s_2^{3}$ large enough,  $m=0,1,2$,
\begin{align}|\delta_m(s)|\leq\Big|\int \chi_1(z,s)\delta(z,s)k_m(z)\mu(z)dz\Big|\leq\int|\delta(z,s)|(1+|z|^2)\mu(z)dz\leq C\frac{s-\sigma}{s^{3+1/2}},\nonumber\end{align}
\begin{align}|\delta_-(s)|\leq| \chi_1(z,s)\delta(z,s)-\sum_{i=0}^2 \delta_i(s)k_i(z)|\leq C\frac{s-\sigma}{s^{3+1/2}}(1+|z|^3)\leq C\frac{s-\sigma}{s^{2+\varepsilon}}(1+|z|^3),\nonumber\end{align}
and
\begin{align}|\delta_e(s)|\leq| (1-\chi_1(z,s))\delta(z,s)|\leq  C\frac{s-\sigma}{s^{1/2+\varepsilon}}.\nonumber\end{align}

Iv) Estimate of $\widetilde{\delta}$.

Since $s\in [\sigma,\sigma+\rho^*)$,  $y(t)\in S^*(t)$, where $t=T-e^{-s}$, we use
 Lemma \ref{BK} and get that \begin{align}&|K(s,\tau)\partial_zG|\nonumber\\\leq & Ce^{s-\tau}\Big\{\frac{\frac{1}{\varepsilon_0}e^{-\frac{(p+1)s}{4(p-1)}}}{\sqrt{1-e^{-(s-\tau)}}}+\frac{s-\tau}{s}(1+s-\tau)
 \Big((1+e^{\frac{s-\tau}{2}})e^{-\frac{s}{2(p-1)}}
+\frac{1}{\varepsilon_0}e^{\frac{(s-\tau)}{2}}e^{-\frac{(p+1)s}{4(p-1)}}\Big)\Big\}.\nonumber\end{align}
Making a change of variables, we write the integral
as
\[
\int_\sigma^s \frac{1}{\sqrt{1-e^{-(s-\tau)}}}d\tau=-\int_{s-\sigma}^0 \frac{1}{\sqrt{1-e^{-\xi}}}d\xi=\int^{s-\sigma}_0 \frac{1}{\sqrt{1-e^{-\xi}}}d\xi=:
\int_0^\eta \frac{1}{\sqrt{1-e^{-\xi}}}d\xi
\]
 with $\xi = s-\tau$.
Since
\[
\frac{1}{\sqrt{1-e^{-\xi}}}\sim \frac 1{\sqrt \xi}\mbox{ as }\xi\to 0
\]
and
$$\lim\limits_{\xi\rightarrow 0}\frac{\frac 1{\sqrt \xi}}{\frac{1}{\sqrt{1-e^{-\xi}}}}=\lim\limits_{\xi\rightarrow 0}\sqrt{\frac{1-e^{-\xi}}{\xi}}=1,$$
there exists $\xi_0>0$ such that
\[
\forall \xi\in (0, \xi_0],\;\; \frac{1}{\sqrt{1-e^{-\xi}}} \le \frac 2{\sqrt \xi}.
\]
Note also that there is $C_0>0$ such that
\[
\forall \xi\ge \xi_0,\;\;  \frac{1}{\sqrt{1-e^{-\xi}}} \le C_0.
\]
Therefore,
\begin{align*}
  &\int_0^\eta \frac{1}{\sqrt{1-e^{-\xi}}}d\xi
  = \int_0^{\min(\eta,\xi_0)} \frac{1}{\sqrt{1-e^{-\xi}}}d\xi
  +\int_{\min(\eta,\xi_0)}^\eta \frac{1}{\sqrt{1-e^{-\xi}}}d\xi\\
  \le &\int_0^{\min(\eta,\xi_0)} \frac 2{\sqrt \xi} d\xi
  +\int_{\min(\eta,\xi_0)}^\eta C_0 d\xi
  =4\sqrt {\min(\eta,\xi_0)}+C_0[\eta-{\min(\eta,\xi_0)}]
  \nonumber\\ \le & 4\sqrt \eta+C_0\eta,
\end{align*}
which implies
\begin{align*}
  &\int_\sigma^{s} \frac{1}{\sqrt{1-e^{-(s-\tau)}}}d\tau
  \le 4\sqrt{s-\sigma} +C_0(s-\sigma).
\end{align*}
Hence,
\begin{align}&\int_\sigma^s|K(s,\tau)\partial_zG(z,\tau)d\tau|\nonumber\\\leq & Ce^{s-\sigma}\Big\{((s-\sigma)+\sqrt{s-\sigma})\frac{1}{\varepsilon_0}e^{-\frac{(p+1)s}{4(p-1)}}+\frac{(s-\sigma)^2}{s}
(1+s-\sigma)(1+e^{\frac{s-\sigma}{2}})e^{-\frac{s}{2(p-1)}}
\nonumber\\&+\frac{(s-\sigma)^2}{\varepsilon_0s}(1+s-\sigma)e^{\frac{s-\sigma}{2}}e^{-\frac{(p+1)s}{4(p-1)}}\Big\}\leq \frac{C((s-\sigma)+\sqrt{s-\sigma})}{s^{3+1/2}},\nonumber\end{align}
if  $s_0\geq s_2^{4}(\varepsilon_0,\rho^*)$ with $s_2^{4}$ large enough.

Then we can get the same estimate on $\widetilde{\delta}$ as $\delta$. This completes the proof of Lemma \ref{main}\ \ \ \#

{\bf Lemma \ref{main} implies Proposition \ref{qa}}
Let $K_0\geq K_{02}$, $\varepsilon_0\in(0,\delta_0]$ and let $\widetilde{A}$ be a positive number. Let $A\geq\widetilde{A}_2(\widetilde{A},K_0,\varepsilon_0)$, where $\widetilde{A}_2(\widetilde{A},K_0,\varepsilon_0)$ will be defined later. Let $s_0\geq\widetilde{s}_2(\widetilde{A},A,K_0,\varepsilon_0)$, where $\widetilde{s}_2(\widetilde{A},A,K_0,\varepsilon_0)$ will be defined later.  Let $q$ be a solution of
(\ref{2.20}) satisfy (\ref{property}), and $s_1\geq s_0$. Assume in addition $\forall s\in[s_0,s_1]$, $y(t)\in S^*(t)$ with $t=T-e^{-s}$.

We want to prove that $\forall s\in[s_0,s_1]$, (\ref{property1}) holds.

Let $\rho_1\geq \rho_2$ be two positive numbers (to be fixed in terms of $A$ later). It is then enough to prove (\ref{property1}), on one hand for $s-s_0\leq \rho_1$,
and on the other hand, for $s-s_0\geq \rho_2$. We suppose $A\geq A_2$, $s_0\geq\max(s_2(A,\widetilde{A},K_0,\varepsilon_0,\rho_1), s_2(A,\widetilde{A},K_0,\varepsilon_0,\rho_2))$.

Case 1: $s-s_0\leq \rho_1$.

Since we have $\forall \tau\in [s_0,s]$, $q(\tau)\in V_{K_0,A}(\tau)$, we apply Lemma \ref{main} (IIi), (Iii), (Iiii)), with $A$, $\rho^*=\rho_1$,
and $\rho=s-s_0$. From (\ref{q}), it holds that
\begin{align}&|q_2(s)|\nonumber\\\leq &\widetilde{A}(\log s_0)s^{-2} + C_1\max\{A,\widetilde{A}\}e^{(s-s_0)}(s-s_0)s^{-3}+(s-s_0)C_1s^{-3}+C_1\frac{s-s_0+\sqrt{s-s_0}}{s^{3+1/2}},\nonumber\\
&|q_-{(z,s)}|\nonumber\\\leq &C_1 \widetilde{A}s^{-2}(1+|z|^3)+(s-s_0)(1+|z|^3){s^{-2-\varepsilon}}+C_1 (s-s_0) s^{-2}(1+|z|^3)\nonumber\\&+C_1\frac{s-s_0+\sqrt{s-s_0}}{s^{2+\varepsilon}}(1+|z|^3),\nonumber\\
&\|q_e(s)\|_{L^\infty}\nonumber\\\leq &C_1(1+e^{(s-s_0)})\widetilde{A}s^{-1/2}+(s-s_0)s^{-\frac{1}{2}-\varepsilon}+(s-s_0)s^{-3/4}+ C_1\frac{s-s_0+\sqrt{s-s_0}}{s^{1/2+\varepsilon}}.\nonumber\end{align}
To have (\ref{property1}), it is enough to satisfy
\begin{align}\label{e2}\widetilde{A}(\log s_0)s^{-2}\leq \frac{A^2}{2}(\log s)s^{-2},\nonumber\\
 C_1 \widetilde{A}s^{-2}+C_1(s-s_0)s^{-2}\leq \frac{A}{4}s^{-2},\nonumber\\
 C_1(1+e^{(s-s_0)})\widetilde{A}s^{-1/2}\leq \frac{A^2}{4}s^{-1/2},
\end{align}
on one hand, and
\begin{align}\label{e3}&C_1\max\{A,\widetilde{A}\}e^{(s-s_0)}(s-s_0)s^{-3}+(s-s_0)C_1s^{-3}+C_1\frac{s-s_0+\sqrt{s-s_0}}{s^{3+1/2}}\nonumber\\\leq &\frac{A^2}{2}(\log s)s^{-2}-s^{-3},\nonumber\\
 &C_1\frac{s-s_0+\sqrt{s-s_0}}{s^{2+\varepsilon}}\leq \frac{A}{4}s^{-2},\nonumber\\
&(s-s_0)s^{-3/4}+ C_1\frac{s-s_0+\sqrt{s-s_0}}{s^{1/2+\varepsilon}}\leq \frac{A^2}{4}s^{-1/2},
\end{align}
on the other hand.

If we restrict $\rho_1$ to satisfy $C_1\rho_1\leq \frac{A}{8}$,  $C_1\widetilde{A}e^{\rho_1}\leq \frac{A^2}{8}$ (which is possible if we fix
$\rho_1=\frac{3}{2}\log A$ for $A$ large), and $A$ to satisfy $\widetilde{A} \leq A$, $\widetilde{A}\leq \frac{A^2}{2}$, $C_1\widetilde{A}\leq \frac{A}{8}$
and $C_1\widetilde{A}\leq \frac{A^2}{8}$, that is $A\geq \widetilde{A}_2^1(\widetilde{A},K_0, \varepsilon_0)$, (note that $C_1$ depends
on $K_0$ and $\varepsilon_0$), then since $s-s_0\leq \rho_1$, (\ref{e2}) holds.

With this value of $\rho_1$, (\ref{e3}) will be satisfied if the following is true:
\begin{align}\label{e4}&C_1AA^{\frac{3}{2}}\frac{3}{2}\log A s^{-3}+\frac{3}{2}\log AC_1s^{-3}+C_1\Big(\frac{3}{2}\log A+\sqrt{\frac{3}{2}\log A}\Big){s^{-(3+1/2)}}\nonumber\\\leq
&\frac{A^2}{2}(\log s)s^{-2}-s^{-3},\nonumber\\
 &C_1\frac{\frac{3}{2}\log A+\sqrt{\frac{3}{2}\log A}}{s^{2+\varepsilon}}\leq \frac{A}{4}s^{-2},\nonumber\\
&(\frac{3}{2}\log A)s^{-3/4}+ C_1\frac{\frac{3}{2}\log A+\sqrt{\frac{3}{2}\log A}}{s^{1/2+\varepsilon}}\leq \frac{A^2}{4}s^{-1/2}.
\end{align}
The first one is equal to
$$C_1A^{\frac{5}{2}}\frac{3}{2}\log A s^{-3}+\frac{3}{2}\log AC_1s^{-3}+C_1\Big(\frac{3}{2}\log A+\sqrt{\frac{3}{2}\log A}\Big){s^{-(3+1/2)}}+s^{-3}\leq \frac{A^2}{2}(\log s)s^{-2},$$
which is possible if $s_0\geq \widehat{s}_2^1(A,\widetilde{A},K_0, \varepsilon_0)$.

Case 2: $s-s_0\geq \rho_2$.

Since we have $\forall \tau\in [\sigma,s]$, $q(\tau)\in V_A(\tau)$, we apply Part I) of Lemma \ref{main} with
$\rho=\rho^*=\rho_2$, $\sigma=s-\rho_2$. From (\ref{q}),
\begin{align}\label{e4}&|q_2(s)|\leq {A}^2(\log (s-\rho_2))s^{-2} + C_2Ae^{\rho_2}\rho_2s^{-3}+C_2\rho_2s^{-3}+C_2\frac{\rho_2+\sqrt{\rho_2}}{s^{3+1/2}},\nonumber\\
&|q_-{(z,s)}|\leq C_2(As^{-2}e^{-\rho_2/2}+A^2e^{-\rho_2^2} s^{-2})(1+|z|^3)\nonumber\\&+C_2(\rho_2+\sqrt{\rho_2})(1+|z|^3){s^{-2-\varepsilon}}+C_2 \rho_2 s^{-2}(1+|z|^3),\nonumber\\
&\|q_e(s)\|_{L^\infty}\leq C_2e^{\rho_2}As^{-1/2}+C_2 A^2s^{-1/2}e^{-\rho_2/p}+\rho_2s^{-3/4}+ C_2\frac{\rho_2+\sqrt{\rho_2}}{s^{1/2+\varepsilon}}.\end{align}

To obtain (\ref{property1}), it is enough to have
\begin{align}f_{A,\rho_2}(s)\geq 0,\nonumber\\C_2(e^{-\rho_2/2}A
+A^2e^{-\rho_2^2}+\rho_2)\leq \frac{A}{4},\nonumber\\C_2(A^2e^{-\rho_2/p}+e^{\rho_2}A)\leq \frac{A^2}{4},
\end{align}
with
$$f_{A,\rho_2}= A^2\log ss^{-2}-s^{-3}-\Big[A^2\log (s-\rho_2)s^{-2}+C_2(Ae^{\rho_2}+1)\rho_2s^{-3}+C_2\frac{\rho_2+\sqrt{\rho_2}}{s^{1/2+3}}\Big],$$
on one hand, and
\begin{align}\label{e5}C_2(\rho_2+\sqrt{\rho_2}){s^{-2-\varepsilon}}\leq \frac{A}{4}s^{-2},\nonumber\\
\rho_2s^{-3/4}+C_2\frac{\rho_2+\sqrt{\rho_2}}{s^{1/2+\varepsilon}}\leq \frac{A^2}{4}s^{-1/2},\end{align}
on the other hand.

Now, it is convenient to fix the value $\rho_2$ such that $C_2Ae^{\rho_2}=\frac{A^2}{8}$, that is $\rho_2=\log \frac{A}{8C_2}$.
The conclusion follows from this choice, for $A$ large. Indeed, for arbitrary $A$, we write
\begin{align}&\Big|f_{A,\log \frac{A}{8C_2}}-s^{-3}\Big[A^2\log \frac{A}{8C_2}-1-C_2\Big(A \frac{A}{8C_2}+1\Big)\log \frac{A}{8C_2}\Big]\Big|\nonumber\\
=&\Big|A^2\log ss^{-2}-s^{-3}-\Big[A^2\log \Big(s-\log \frac{A}{8C_2}\Big)s^{-2}+C_2\Big(A\frac{A}{8C_2}+1\Big)\log \frac{A}{8C_2}s^{-3}\nonumber\\&+C_2\frac{\log \frac{A}{8C_2}+\sqrt{\log \frac{A}{8C_2}}}{s^{1/2+3}}\Big]-s^{-3}\Big[A^2\log \frac{A}{8C_2}-1-C_2\Big(A \frac{A}{8C_2}+1\Big)\log \frac{A}{8C_2}\Big]\Big| \nonumber\\
=&\Big|A^2\log ss^{-2}-A^2\log \Big(s-\log \frac{A}{8C_2}\Big)s^{-2}-C_2\frac{\log \frac{A}{8C_2}+\sqrt{\log \frac{A}{8C_2}}}{s^{1/2+3}}-s^{-3}A^2\log \frac{A}{8C_2}\Big|. \nonumber
\end{align}
Note that when $|x|\leq \widetilde{\varepsilon}$ with $\widetilde{\varepsilon}$  small enough, there exists a constant $C_3>0$
such that $|\log(1+x)-x|\leq C_3x^2$. Thus, we have

\begin{align}
&\Big|A^2\log ss^{-2}-A^2\log \Big(s-\log \frac{A}{8C_2}\Big)s^{-2}-C_2\frac{\log \frac{A}{8C_2}+\sqrt{\log \frac{A}{8C_2}}}{s^{1/2+3}}-s^{-3}A^2\log \frac{A}{8C_2}\Big| \nonumber\\
\leq &A^2s^{-2}\Big|\log s-\log \Big(s-\log \frac{A}{8C_2}\Big)-\log \frac{A}{8C_2}s^{-1}\Big|+C_2\frac{\log \frac{A}{8C_2}+\sqrt{\log \frac{A}{8C_2}}}{s^{1/2+3}}\nonumber\\
=&A^2s^{-2}\Big|-\log \Big(1-\frac{1}{s}\log \frac{A}{8C_2}\Big)-\log \frac{A}{8C_2}s^{-1}\Big|+C_2\frac{\log \frac{A}{8C_2}+\sqrt{\log \frac{A}{8C_2}}}{s^{1/2+3}}\nonumber\\
\leq &\frac{A^2}{s^2}C_3\frac{1}{s^2}\Big(\log \frac{A}{8C_2}\Big)^2+C_2\frac{\log \frac{A}{8C_2}+\sqrt{\log \frac{A}{8C_2}}}{s^{1/2+3}}
\nonumber\\\leq & C_2A^2\frac{\Big(\log \frac{A}{8C_2}\Big)^2}{s^{1/2+3}},\nonumber\end{align}
when $A\geq \widetilde{A}_2^2(K_0,\varepsilon_0)$, and then take $s_0\geq \widehat{s}_2^2(A,K_0,\varepsilon_0)$.

Then we take $A\geq \widetilde{A}_2^3(K_0,\varepsilon_0)$ such that
\begin{align}A^2\log \frac{A}{8C_2}-1-C_2\Big(A \frac{A}{8C_2}+1\Big)\log \frac{A}{8C_2}\geq 1,\nonumber\\C_2\Big(\Big(\frac{A}{8C_2}\Big)^{-\frac{1}{2}}A
+e^{-(\log \frac{A}{8C_2})^2}A^2+\Big(\log \frac{A}{8C_2}\Big)\Big)\leq \frac{A}{4}\nonumber,\end{align}\begin{align}C_2\Big(A^2(\frac{A}{8C_2}\Big)^{-\frac{1}{p}}+\frac{A}{8C_2}A\Big)\leq \frac{A^2}{4}.\nonumber\end{align}

After, we introduce $s\geq s_0\geq \widehat{s}_2^3(A,K_0,\varepsilon_0)$, we have
$$C_2A^2\frac{(\log \frac{A}{8C_2})^2}{s^{1/2+3}}\leq \frac{1}{2}s^{-3}$$
and (\ref{e5}) holds.

This way, for $A\geq \widetilde{A}_2(\widetilde{A},K_0,\varepsilon_0)$ and $s\geq s_0\geq \widehat{s}_2(\widetilde{A},A,K_0,\varepsilon_0)$, which concludes Case 2, and we complete the
proof of Proposition \ref{qa}.\ \ \ \#

{\bf Part 2: Estimate in $\mathcal{R}_2$}

The aim of this part is to show that if
\begin{align}\label{P2}\frac{\epsilon_0}{2}\leq |x-a|,\ x\in \Omega,\ \mbox{then}\ |y(x,t_*)|\leq \frac{\eta_0}{2}\end{align}
provided the parameters satisfy some conditions.

Step 1. Improved estimates in the intermediate region.

Here, we refine the estimates on the solution in the region
\begin{align} \label{bu region} K_0\sqrt{(T-t)|\log(T-t)}|\leq |x-a|\leq \frac{\varepsilon_0}{2}.\end{align}
Given a small $x$, define $t=t_0(x)$ by
\begin{align} \label{bu region1} |x|=K_0\sqrt{(T-t_0(x))|\log(T-t_0(x))}|\end{align}
to see that the solution is in fact flat at that time. Then,  we see that the solution of (\ref{1.6}) remains
flat for later time. More precisely, we claim the following:

\begin{lemma}\label{estimate}Assume that $y(t)\in S^*(t)$, $t\in[0,t_*]$. Then, there exist $0<\varsigma_0<1$, $K_0^2>0$ and $\widehat{\delta}\in
(0,\delta_0)$ such that for all $K_0\geq K_0^2$, $\varepsilon_0\in (0,\widehat{\delta}]$ and
 $A\geq 1$, there exists $\widehat{s}_3(K_0,\varepsilon_0,A)$
such that if $s_0\geq \widehat{s}_3$ and $0<\eta_0\leq 1$, then for  $x_0\in  \{x_0\in \mathbb{R}:\ 0<|x_0-a|<\varepsilon_0\leq \widehat{\delta}\}$
  and $x=x_0+\xi\sqrt{T-t_0(x_0)}$ with $|\xi|\leq |\log(T-t_0(x_0))|^{1/4}$, it holds that
$$\forall t_0(x_0)\leq  t\leq t_*,\ \Big|\frac{y(x,t)}{y^*(x_0)}-\frac{U_{K_0}(\tau)}{U_{K_0}(1)}\Big|\leq \frac{C}{|\log|x_0-a||^{\varsigma_0}},$$
where
\begin{align}y^*(x_0)=\Big[\frac{(p-1)^2|x_0-a|^2}{8p|\log|x_0-a||}\Big]^{-\frac{1}{p-1}},\end{align}
and
\begin{align}U_{K_0}(\tau)=\kappa((1-\tau)+\frac{p-1}{4p}K_0^2)^{-\frac{1}{p-1}}.\end{align}
Moreover, $|y(x_0,t)|\leq C_4(K_0)|y^*(x_0)|$, $t\in[0,t_*]$, where $C_4(K_0)$ is a constant depending on $K_0$.
\end{lemma}

 The proof of Lemma \ref{estimate} is a small adaptation of the treatment of estimate (1.11)  in \cite{Abdelhedi}, and that for details, one may see \cite{Abdelhedi}.

Step 2. A parabolic estimate in Region $\mathcal{R}_2$.

We recall from the definition of $S^*(t)$ that
$\mbox{for all}\ \frac{\varepsilon_0}{2}\leq |x-a|,\ x\in \Omega,\ |y(x,t)|\leq \eta_0$.

Now we will obtain a parabolic estimate on the solution in $\mathcal{R}_2$.

\begin{proposition}\label{P4}(A parabolic estimate in Region $\mathcal{R}_2$)

For any $\varepsilon\in(0,1]$, any $\varepsilon_0\in(0,\widehat{\delta}]$ and any $\sigma_1\geq 0$, there exist
 a positive constant $C_5$ independent of $\varepsilon,\varepsilon_0,\sigma_1$ and $T_1(\varepsilon,\varepsilon_0,\sigma_1)=\min\Big\{1,\Big(\frac{\varepsilon/2}{C_5e^{(\sigma_1+1)^{p-1}}(\sigma_1+1)
(\frac{1}{\varepsilon_0}+\frac{1}{\varepsilon_0^2})}\Big)^2\Big\}$
such that for all $\overline{t}$ with $0<\overline{t}\leq T_1$, if $y$ is a solution of
$$y_t=y_{xx}+\chi_\omega|y|^{p-1}y,\ \forall \ x\in \Omega,\ t\in [0,\overline{t}]$$
which satisfies

(i) for $$|x-a|\in \big[\frac{\varepsilon_0}{4},\frac{\varepsilon_0}{2}\big],\ |y(x,t)|\leq \sigma_1,$$

(ii) for $$|x-a|\geq\frac{\varepsilon_0}{4},\ x\in \Omega,\ y(x,0)=0.$$
Then for any $t\in[0,\overline{t}]$ and any $x\in \Omega$ with $\frac{\varepsilon_0}{2}\leq |x-a|,$ it holds that
$$|y(x,t)|< \varepsilon.$$
\end{proposition}

Proof. Recall that $\overline{y}$ defined in
\begin{align}\overline{y}(x,t)=& y(x,t)\overline{\chi}(x),\ x\in\Omega,\end{align}
where $\overline{\chi}(x)=1-\chi_0(\frac{4(x-a)}{\varepsilon_0})$, which satisfies
\begin{align}\label{1.8}
\partial_t \overline{y}=\partial_{xx}\overline{y}+\chi_\omega|y|^{p-1}\overline{y}-2\partial_x(\overline{\chi}' y)+\overline{\chi}{'}{'}y\nonumber.
\end{align}
Therefore, by (ii), we can write
\begin{align}
\|\overline{y}(t)\|_{L^\infty(\Omega)}\leq \int_0^t S(t-t')\big[\chi_\omega|y|^{p-1}I_{|x-a|\geq\frac{\varepsilon_0}{4}}\overline{y}\nonumber\\
-2\partial_x(\overline{\chi}' I_{|x-a|\geq\frac{\varepsilon_0}{4}}y)+\overline{\chi}{'}{'}I_{|x-a|\geq\frac{\varepsilon_0}{4}}y\big]dt',
\end{align}
where $S(t)$ is the heat kernel.

We proceed by contradiction, and assume that there exists $\overline{t}$ with $0<\overline{t}\leq T_1$ with $T_1(\varepsilon,\varepsilon_0,\sigma_1)=\min\Big\{1,\Big(\frac{\varepsilon/2}{C_5e^{(\sigma_1+1)^{p-1}}(\sigma_1+1)
(\frac{1}{\varepsilon_0}+\frac{1}{\varepsilon_0^2})}\Big)^2\Big\}$, such that the conclusion does not hold for all $t\in[0,\overline{t}]$, where
$C_5$ will be defined later.

Using (ii) and the continuity of $y$,
this means that there is $\widehat{t}\in (0,\bar t]$ such that the
conclusion holds for all $t\in [0, \widehat{t})$ and fails at $t=\widehat{t}$. This
means that
\begin{equation}\label{five}
\|\overline{y}(\widehat{t})\|_{L^\infty(\frac {\epsilon_0}2 \le |x-a|,\ x\in \Omega)}=\varepsilon.
\end{equation}
Therefore, since $\overline{\chi}'$ and $\overline{\chi}''$ are supported by $\{\frac{\varepsilon_0}{4}\leq |x-a|\leq\frac{\varepsilon_0}{2}\}$
and satisfy $|\overline{\chi}'|\leq C/\varepsilon_0$, $|\overline{\chi}''|\leq C/\varepsilon_0^2$, it holds that
for  $(x,t)\in\Omega\times[0,\widehat{t}]$,
\begin{align}
&\|\overline{y}(t)\|_{L^\infty(\Omega)}\nonumber\\\leq &(\sigma_1+1)^{p-1} \int_0^t \|\overline{y}(t')\|dt'+\frac{C(\sigma_1+1)}{\varepsilon_0}\int_0^t\frac{dt'}{\sqrt{t-t'}}+\frac{C(\sigma_1+1)}{\varepsilon_0^2}\int_0^t dt'\nonumber\\
\leq &(\sigma_1+1)^{p-1} \int_0^t \|\overline{y}(t')\|dt'+\frac{C(\sigma_1+1)}{\varepsilon_0}\sqrt{\widehat{t}}+\frac{C(\sigma_1+1)}{\varepsilon_0^2}\widehat{t},\ t\in[0,\widehat{t}].\nonumber
\end{align}
In the above estimate, when it comes to bounding $|y|1_{\{|x-a|\ge
  \frac{\epsilon_0}4\}}$, we will bound it by $\sigma_1$ if
$|x-a|\le \frac{\epsilon_0}2$, and by $\varepsilon$ if $|x-a|\ge
\frac{\epsilon_0}2$. Then by Gronwall estimate, we have
\begin{align}
\|\overline{y}(t)\|_{L^\infty(\Omega)}&\leq Ce^{(\sigma_1+1)^{p-1}}\Big(\frac{C(\sigma_1+1)}{\varepsilon_0}\sqrt{\widehat{t}}+\frac{C(\sigma_1+1)}{\varepsilon_0^2}\widehat{t}\Big)\nonumber\\
&\leq C_5e^{(\sigma_1+1)^{p-1}}(\sigma_1+1)\sqrt{\widehat{t}}(\frac{1}{\varepsilon_0}+\frac{1}{\varepsilon_0^2}),\ t\in[0,\widehat{t}].
\end{align}
Hence, if $T_1(\varepsilon,\varepsilon_0,\sigma_1)=\min\Big\{1,\Big(\frac{\varepsilon/2}{C_5e^{(\sigma_1+1)^{p-1}}(\sigma_1+1)
(\frac{1}{\varepsilon_0}+\frac{1}{\varepsilon_0^2})}\Big)^2\Big\}$  and $\widehat{t}\leq \overline{t}\leq T_1$, we have $\|\overline{y}(t)\|_{L^\infty(\Omega)}\leq \varepsilon/2$,  $t\in[0,\widehat{t}]$.

This contradicts with \eqref{five} and completes the proof of Proposition \ref{P4}.\ \ \ \#

Step 3. Proof of the improvement in (\ref{P2}).

Here we use step 1 and step 2 to prove (\ref{P2}) for a suitable choice of parameters.
 We consider $K_0\geq K_0^2$,
$\varepsilon_0\in
(0,\widehat{\delta}]$, $A\geq 1$, $0<\eta_0\leq 1$, and
$$s_0\geq \widehat{s}_{4}=:\max{\Big\{\widehat{s}_{3}(K_0,\varepsilon_0,A),\ s_{0,1}(K_0,\varepsilon_0,A),
\ -\log\Big(T_1\Big(\frac{\eta_0}{2},\varepsilon_0,C_4(K_0)\Big|y^*(\frac{\varepsilon_0}{4})\Big|\Big)\Big\}},$$
where the different constants are defined in Lemma \ref{initial data1}, Lemma \ref{estimate} and Proposition \ref{P4}. Applying Lemma \ref{estimate}, we see that if $y(t)\in S^*(t)$, $t\in[0,t_*]$,
 $$\forall \frac{\varepsilon_0}{4}\leq|x-a|\leq \frac{\varepsilon_0}{2},\ \forall \ t\in[0,t^*],\ |y(x,t)|\leq C_4(K_0)|y^*({x})|\leq C_4(K_0)|y^*({\frac{\varepsilon_0}{4}})|.$$

Using the proof of Lemma \ref{initial data1}, $\forall \frac{\varepsilon_0}{4}\leq|x-a|$, $x\in \Omega$, $y(x,0)=0$.

Therefore, Proposition \ref{P4} applies with $\varepsilon=\frac{\eta_0}{2}$ and $\sigma_1=C_4(K_0)|y^*({\frac{\varepsilon_0}{4}})|$,
 \begin{align}\label{P5}\forall \frac{\varepsilon_0}{2}\leq|x-a|,\ x\in \Omega, \forall \ t\in[0,t_*],\ |y(x,t)|\leq \frac{\eta_0}{2}, \ t\in[0,t_*].\end{align}

Hence, (\ref{P2}) holds.\\

{\bf 2.3.2 Transverse crossing on $V_{K_0,A}$}

Similarly to Lemma 3.8 in \cite{Zaag1}, we have the following lemma on transverse crossing on $V_{K_0,A}$.

\begin{lemma}\label{lemma 3.2} There exist $A_3>0$ and $K_0^3>0$ such that for any $A\geq A_3\geq 1$,  $K_0\geq K_0^3$,
 $\varepsilon_0\in(0,\widehat{\delta}]$ and $\eta_0\in(0,1]$, there exists $\widehat{s}_5(K_0,\varepsilon_0,A)$ such that for any $s_0\geq \widehat{s}_5(K_0,\varepsilon_0,A)$,  we have the following properties:
Assume there exists $\widetilde{s}_*\geq s_0$ such that $y(\widetilde{t}_*)\in S(\widetilde{t}_*)$ with $\widetilde{t}_*=T-e^{-\widetilde{s}_*}$ and $(q_0,q_1)(\widetilde{s}_*)\in \partial [-\frac{A}{\widetilde{s}_*^2}, \frac{A}{\widetilde{s}_*^2}]^2$, then there exists $\delta_1>0$
such that  $\forall \delta\in (0,\delta_1)$, $(q_0,q_1)(\widetilde{s}_*+\delta)\not\in [-\frac{A}{(\widetilde{s}_*+\delta)^2}, \frac{A}{(\widetilde{s}_*+\delta)^2}]^2$.
\end{lemma}

\subsection{Topological argument}
In this section, we shall conclude the proof of Theorem \ref{Pros1} by Topological argument.

  We fix $S_0\geq \max\{\widetilde{s}_0, s_{0,1}, \widehat{s}_1,\widehat{s}_4, \widehat{s}_5\}$, $A\geq \max\{A_1,A_3\}$,
  $K_0\geq \max\{K_0^1, K_0^2, K_0^3\}$, $\varepsilon_0\in(0,\widehat{\delta}]$ and $\eta_0\in(0,1]$,
  take $s_0\geq S_0$. We argue by contradiction: According to Lemma \ref{initial data1},  for all $(d_0,d_1)\in D_T$,
  $y_0(\cdot,d_0,d_1)\in S^*(0)$. We suppose then that for each
$(d_0,d_1)\in D_T$, there exists $s>s_0$ such that $y(t) \not\in S(t)$ ($t=T-e^{-s}$). Let $s_*(d_0,d_1)$ be the infinimum of all these $s$.

 Applying Proposition \ref{qb} and
 (\ref{P2}), we see that $y(t^*)$ can leave $S(t^*)$ only by its first two components, hence, $s_*=-\log(T-t^*)$, and
 $$(q_0,q_1)(d_0,d_1,s_*(d_0,d_1))\in \partial S(t^*)(s_*(d_0,d_1)).$$
We see from Definition \ref{s} of $S^*(t)$ that only  the components $q_0(s^*)$ or $q_1(s^*)$ may touch
the boundary of $[-\frac{A}{s_*^2},\frac{A}{s_*^2}]$.

Then we may define the rescaled flow:
$$\Phi: D_T\rightarrow \partial ([-1,1])^2,$$
$$(d_0,d_1)\rightarrow \frac{s_*(d_0,d_1)^2}{A}(q_0,q_1)(d_0,d_1,s_*(d_0,d_1)).$$
In particular, either
\begin{align}\label{Transverse}wq_0(s_*)=\frac{A}{s_*^2},\ \mbox{or}\ wq_1(s_*)=\frac{A}{s_*^2},\end{align}
$w\in\{-1,1\}$, both depending on $(d_0,d_1)$.\\

 Now we claim that

\begin{proposition}
1) $\Phi$ is continuous mapping from $D_T$ to $\partial ([-1,1]^2)$.

2) The restriction of $\Phi$ to $\partial D_T$ is homeomorphic to identity.
\end{proposition}

Proof. By Lemma \ref{lemma 3.2}
and Lemma \ref{initial data1}, we can use the similar techniques in the proof of Proposition 3.10 in \cite{Zaag1} to prove this proposition.

Form that, a contradiction follows by Index Theory.  This means that there exists $(d_0,d_1)\in D_T$ such that $y(t)\in S^*(t)$
and $T=t_*(d_0,d_1)$.

Let us fix $K_0>0$, $\varepsilon_0>0$, $A>0$, $0<\eta_0\leq 1$ and $T>0$ so that  all the statements of Section 2.1-2.3 apply.  Hence, for some $(d_0,d_1)\in D_T$, equation (\ref{1.6})
with initial data (\ref{initial data}) has a solution $y$ such that $T=t_*(d_0,d_1)$ and for any $t\in[0,T)$, $y(t)\in S^*(K_0,\varepsilon_0,A_0,\eta_0,T,t)$.
By  (\ref{W4}), we see that for any $s\geq-\log T$ and for any $z\in \mathbb{R}$,
$$|q(z,s)|\leq \frac{CA^2}{\sqrt{s}}.$$
By (\ref{w}), (\ref{w1}) and (\ref{w2}), it holds that  for any $s\geq-\log T$ and for any $z$ with $|z|\leq \varepsilon_0e^{s/2},$
$$\Big|W(z,s)-f(\frac{z}{\sqrt{s}})\Big|\leq \frac{CA^2}{\sqrt{s}}+\frac{C}{s}.$$
Hence, if $T$ is small enough, we have that for any
$t\in[0,T)$ and for any $x\in \Omega$ with $|x-a|\leq \varepsilon_0 $,
\begin{align}\label{e*}\Big|(T-t)^{1/(p-1)} y(x,t)-f\Big(\frac{(x-a)(T-t)^{-1/2}}{\sqrt{|\log(T-t)|}}\Big)\Big|\leq \frac{C(A)}{\sqrt{|\log(T-t)|}}.\end{align}
This implies (\ref{1.3}).

Now just take $a_j \equiv a$ and $t_j = T-\epsilon_j$, $j=1,2,\cdots$, with $\epsilon_j$
being any sequence converging to $0$.
Then, \begin{align}\nonumber\Big| \frac{y(a,t)}{\kappa(T-t_j)^{-1/(p-1)}}-1\Big|\leq \frac{C(A)}{\sqrt{|\log(T-t_j)|}}\rightarrow 0,\ \mbox{as}\ j\rightarrow +\infty,\end{align}
that is, $|y(a,t_j)|\sim \kappa(T-t_j)^{-1/(p-1)}$ as $j\rightarrow +\infty$. Hence, $y$ blows up at point $a$ in time $T$.

It remains to prove that for any $x\neq a$ is not a blowup point.

We know from (ii) in Definition \ref{s} that for all $x\in \Omega$ with $|x-a|\geq \frac{\varepsilon_0}{2},\ |y(x,t)|\leq \eta_0$.
Thus,  any $x_0\in \Omega$ with $|x_0-a|\geq \frac{\varepsilon_0}{2}$ is not a blowup point. Now, if $0<|x_0-a|\leq \varepsilon_0/2$, the following result from Giga and Kohn \cite{Giga-1} allows us to conclude
Theorem \ref{Pros1}.

\begin{proposition}\label{GK}
For all $C_0>0$, there is $\eta_0>0$ such that if $v(\xi,\tau)$ solves
$$|v_t-\Delta v|\leq C_0(1+|v|^p),$$
and satisfies
$$|v(\xi,\tau)|\leq \eta_0(T-t)^{-\frac{1}{p-1}},$$
for all $(\xi,\tau)\in B(a,r)\times [T-r^2,T)$ for some $a\in \mathbb{R}$ and $0<r\leq 1$, then $v$ does not blows up at $(a,T)$.
\end{proposition}

  Indeed, since $0<|x_0-a|\leq \varepsilon_0/2$, it follows from (\ref{e*}),
\begin{align}\sup\limits_{|x-x_0|\leq\frac{|x_0-a|}{2}}|(T-t)^{1/(p-1)} y(x,t)\Big|\leq \Big|f\Big(\frac{|x_0-a|/2(T-t)^{-1/2}}{\sqrt{|\log(T-t)|}}\Big)\Big|+ \frac{C(A)}{\sqrt{|\log(T-t)|}}\rightarrow 0,\nonumber\end{align}
as $t\rightarrow T$.
Therefore, applying Proposition \ref{GK}, we see that $x_0$ is not a blowup point.

This completes the proof of Theorem \ref{Pros1}.\ \ \ \#

\section{Proof of Theorem \ref{Main Tho}}

Given $p>1$. By Theorem \ref{Pros1}, we  conclude  that for each $a\in \omega$, there exists $T_0>0$ such that for each $\widetilde{T}\in (0,T_0)$,
there exists an initial data $\widetilde{y}_0\in C_0^\infty(\omega)$ such that the  solution
$y$ to  (\ref{1.6}) blows up in time $\widetilde{T}$ and   has unique blowup point $a$.

Then,  given  $(a,T)\in \omega\times(0,+\infty)$,  we take a time  $T_1\in(0, \min\{T_0,T/2\})$. Then there exists an initial data $\widetilde{y}_0\in C_0^\infty(\omega)$ such that   solution
$y$ to  (\ref{1.6})  has unique blowup point $(a,T_1)$. We set $$\overline{y}(x,t):=y(x, t-T+T_1), \ (x,t)\in \Omega
 \times [T-T_1, T).$$
 Then, $\overline{y}$ satisfies the following system,
\begin{eqnarray}\nonumber
\left\{\begin{array}{ll} \overline{y}_t-
\Delta \overline{y}=\chi_\omega |\overline{y}|^{p-1}\overline{y},&x\in \Omega,\ t\in(T-T_1, T),\\
\overline{y}=0,& x\in \partial\Omega,\  t\in(T-T_1,T),\\
\displaystyle  \overline{y}(x,T-T_1)=\widetilde{y}_0(x),&
x\in \Omega.
\end{array}\right.
\end{eqnarray}
Moreover, $\overline{y}$ has unique blowup point $(a,T)$.

On the other hand, for each $y_0\in H_0^1(\Omega)$,  take
the following system in consideration,
\begin{eqnarray}\label{13}
\left\{\begin{array}{ll} z_t-
\Delta z=\chi_\omega v,&(x,t)\in \Omega\times  (0,T-T_1),\\
z=0,&(x,t)\in \partial\Omega\times  (0,T-T_1),\\
\displaystyle  z(x,0)=y_0(x)-\widetilde{y}_0(x),&x\in \Omega,
\end{array}\right.
\end{eqnarray}
where $v\in L^2(0,T-T_1;H)$. Since $y_0-\widetilde{y}_0\in H_0^1(\Omega)$, it is well known that for each $v\in L^2(0,T-T_1;H)$, there exists a unique  solution $z$ in $C([0,T-T_1];H_0^1(\Omega))$ to (\ref{13}).
(\ref{13}) can be equivalently written as
\begin{eqnarray}\label{13-1}
\left\{\begin{array}{ll} z'(t)=Az(t)+Bv(t), \ t\in (0,T-T_1)
,\\
z(0)=y_0-\widetilde{y}_0(x).
\end{array}\right.
\end{eqnarray}

Consider the following optimal control problem,
$$(\mathcal{P})\ \min\Big\{\int_0^{T-T_1}\|v(t)\|_{H}^2dt;\ z'=Az+Bv,\ z(0)=y_0-\widetilde{y}_0,\ z(T-T_1)=0\Big\},$$
and the  Riccati system (\ref{13-2}).

By Theorem 2.1 in \cite{Sirbu} and its proof, there exists a unique mild solution $P\in C_S([0,T-T_1); \Sigma^+(H))$ to problem
(\ref{13-2}). Moreover, $P$ satisfies  $\lim\limits_{t\rightarrow T-T_1}\langle P(t)z(t), z(t) \rangle=0$, for every mild solution $z$ of the
state system $z'=Az+Bv$, $z(t_0)=z_0$ with $0\leq t_0<T-T_1$, $z(T-T_1)=0$ and $v\in L^2(t_0,T-T_1;H)$.

Moreover, $v(t)=-B^*P(t)z(t)$, $t\in [0,T-T_1)$,  is the optimal feedback control for problem ($\mathcal{P}$).

Set $\widehat{y}(t):=z(t)+\widetilde{y}_0(x)$, $t\in [0,T-T_1)$ and $$u_1(t):=v(t)-\Delta \widetilde{y}_0(x)=-B^*P(t)(\widehat{y}(t)-\widetilde{y}_0(x))-
\Delta \widetilde{y}_0(x),\ t\in[0,T-T_1).$$  It
holds by the construction of $\widetilde{y}_0(x)$ that
\begin{eqnarray}\nonumber\Delta \widetilde{y}_0(x)=\chi_\omega \Delta \widetilde{y}_0(x).\end{eqnarray}
Then,  it follows that $\widehat{y}$ is the solution to the following system,
\begin{eqnarray}\nonumber
\left\{\begin{array}{ll} \widehat{y}_t-
\Delta \widehat{y}=\chi_\omega u_1,&x\in \Omega\times  (0,T-T_1),\\
\widehat{y}=0,& x\in \partial\Omega\times  (0,T-T_1),\\
\displaystyle  \widehat{y}(x,0)=y_0(x),&
x\in \Omega,
\end{array}\right.
\end{eqnarray}
and \begin{eqnarray}\nonumber\widehat{y}(T-T_1)=\widetilde{y}_0(x).\end{eqnarray}

Set
$$y(x,t):=\left\{\begin{array}{ll}\widehat{y}(x,t),\ &(x,t)\in \Omega\times  [0,T-T_1), \\
\overline{y}(x,t), \ &(x,t)\in \Omega\times [T-T_1,T),
\end{array}\right.$$
and
$$u(x,t):=\left\{\begin{array}{ll}-B^*P(t)({y}(t)-\widetilde{y}_0)(x)-\Delta \widetilde{y}_0(x),\ &(x,t)\in \Omega\times  (0,T-T_1), \\
|y|^{p-1}y(x,t), \ &(x,t)\in \Omega\times [T-T_1,T).
\end{array}\right.$$

Then, $y$ is the solution to (\ref{e1.1})
with  the   feedback control $u$, and  it follows that
 $y$ blows up in $T$ and has unique blowup point $a$.

This completes the proof of Theorem \ref{Main Tho}.\ \ \ \#
\section{Proof of Noncontrollability}
Proof of Theorem \ref{Tho 2}. Let $y_0\in H_0^1(\Omega)\bigcap L^\infty(\Omega)$. Suppose that $y$ is a corresponding solution to system (\ref{e1.1}) for some feedback control and belongs to   $C([0,t_{max});L^\infty(\Omega))$, and suppose that   $a\in \Omega\setminus \overline{\omega}$ is the unique blowup point of $y$. Then we can find three balls $B_0$, $B_1$ and $B_2$
with $a\in B_0\subset\subset B_1\subset\subset B_2
\subset\subset  \Omega\setminus \overline{\omega}$. Set $\chi_2\in C^\infty (\Omega\setminus \overline{\omega})$,
\begin{eqnarray}\nonumber\chi_2:=\left\{\begin{array}{ll}1,\ &x\in B_0, \\
0, \ &x\in B_2\setminus \overline{B}_1.
\end{array}\right.\end{eqnarray}
Let $\varphi=\chi_2 y$. We get \begin{eqnarray}
\left\{\begin{array}{ll} \varphi_t-
\Delta \varphi=-2\nabla \chi_2\nabla y-\Delta\chi_2  y=-2\nabla(y\cdot\nabla \chi_2)+\Delta\chi_2  y,\ &x\in B_2,\ \ \ t\in[0,t_{max}),\\
\varphi=0,& x\in \partial B_2,\  t\in [0,t_{max}),\\
\displaystyle  \varphi(x,0)=(\chi_2y_0)(x),&
x\in B_2.
\end{array}\right.
\end{eqnarray}
Then semigroup representation formula for $\varphi$ gives
\begin{eqnarray}\varphi(t)=e^{t}(\chi_2y_0)+\int_0^t e^{(t-s)\Delta}(-2\nabla\cdot(y\nabla \chi_2)+\Delta\chi_2  y)ds,\ t\in[0,t_{max}).\end{eqnarray}
Then, we have
\begin{align}\label{non}\|\varphi(t)\|_{L^\infty}\leq &C\|y_0\|_{L^\infty}+C\int_0^t\big( (t-s)^{-1/2}\|y\nabla \chi_2\|_{L^\infty}
\nonumber\\&+\|\Delta\chi_2  y\|_{L^\infty}\big)ds,\ t\in[0,t_{max}).\end{align}

Since $a\in B_0\subset\subset B_1\subset\subset B_2
\subset\subset  \Omega\setminus \overline{\omega}$ and $a$ is the unique blowup point  of $y$, we have that there exists a ball $B_3$ with $a\in{B}_3\subset\subset B_0$ and a constant
$\overline{C}>0$ such that for any
 $t\in[0,t_{max})$, $\|y(t)\|_{L^\infty}\leq \overline{C}$ in $B_2\setminus \overline{B}_3$. On the other hand, $\nabla \chi_2=0$ and $\Delta\chi_2=0$ in $B_0$. Then by (\ref{non}), it holds that for any
 $t\in[0,t_{max})$,  $\|\varphi\|_{L^\infty}\leq \overline{C}_1$ in $B_2$ for some $\overline{C}_1>0$. This contradicts with
the assumption $a$ is a blowup point.

This comletes the proof of Theorem 1.2.\ \ \ \#
\section{Appendix}

{\bf  Proof of Lemma \ref{BK}.}

 Step 1. Perturbation formula for $K(s,\sigma,z,x)$

Since $\mathcal{L}$ is conjugated to the harmonic oscillator $e^{-\frac{x^2}{8}}\mathcal{L}e^{\frac{x^2}{8}}=\partial^2-\frac{x^2}{16}+1/4+1$, we use the definition of $K$
 and give a Feynman-Kac representation for $K$,
 \begin{align}K(s,\sigma,z,x)=e^{(s-\sigma)\mathcal{L}}(z,x)\int d\mu_{zx}^{s-\sigma}(w)e^{\int_0^{s-\sigma}V(w(\tau),\sigma+\tau)}d\tau,\end{align}
 where $d\mu_{zx}^{s-\sigma}$ is the oscillator measure on the continuous paths $w:[0,s-\sigma]\rightarrow \mathbb{R}$ with $w(0)=x$, $w(s-\sigma)=z$, i.e. the
 Gaussian probability measure with covariance kernel $\Gamma(\tau,\tau')$,
  \begin{align}\label{Gamma}=w_0(\tau)w_0(\tau')+2(e^{-1/2|\tau-\tau'|}-e^{-1/2|\tau+\tau'|}+e^{-1/2|2(s-\sigma)-\tau'
  +\tau|}-e^{-1/2|2(s-\sigma)-\tau'-\tau|}),\end{align}
  which yields $\int d\mu_{zx}^{s-\sigma}w(\tau)=w_0(\tau)$ with
  $$w_0(\tau)=(\sinh \frac{s-\sigma}{2})^{-1}(z\sinh \frac{\tau}{2}+x\sinh\frac{s-\sigma-\tau}{2}).$$
  We have in addition
  $$e^{\theta\mathcal{L}}(z,x)=\frac{e^\theta}{\sqrt{4\pi(1-e^{-\theta})}}\exp\Big[-\frac{(ze^{-\theta/2}-x)^2}{4(1-e^{-\theta})}\Big].$$

  We write from now on $(\psi,\varphi)$ for $\int d\mu \psi(z)\varphi(z)$.
  \begin{lemma} \label{B-K1}$\forall s\geq \sigma\geq 1$ with $s\leq 2\sigma$, the kernel $K(s,\sigma,z,x)$ satisfies
  $$K(s,\sigma,z,x)=e^{(s-\sigma)\mathcal{L}}(z,x)\Big(1+\frac{1}{s}P_1(s,\sigma,z,x)+P_2(s,\sigma,z,x)\Big)$$
  \end{lemma}
  where $P_1$ is a polynomial
  $$P_1(s,\sigma,z,x)=\sum_{m,n\geq 0,m+n\leq2}p_{m,n}(s,\sigma)z^mx^n$$
  with $|p_{m,n}(s,\sigma)|\leq C(s-\sigma)$ and
  $$|P_2(s,\sigma,z,x)|\leq C(s-\sigma)(1+s-\sigma)s^{-2}(1+|z|+|x|)^4.$$
  Moreover, $|(k_2,(K(s,\sigma)-(\sigma s^{-1})^2)h_2|\leq C(s-\sigma)(1+s-\sigma)s^{-2}.$
 See Lemma 5 in \cite{Bricmont}.

  Step 2. Conclusion of the proof of Lemma \ref{BK}.

  a), b) and d) see Lemma 3.13 in \cite{Zaag1}, Lemma 3.11  in \cite{Zaag2} and Lemma 2 in \cite{Bricmont}, respectively.

  Proof of c). We consider $K_0>0$, $A'>0$, $A^{''}>0$, $A^{'''}>0$ and $\rho*>0$. Let $s_0\geq \rho*$, $\sigma\geq s_0$ and $q(\sigma)$ satisfying
  (\ref{qsigma}).
    We estimate $$\alpha(s)=K(s,\sigma)q(\sigma)=\sum_{m=0}^{2}\alpha_m(s) h_m(z)+\alpha_{-}(z,s)+\alpha_e(z,s),$$
    for each $s\in [\sigma,\sigma+\rho^*]$.

    Since $\sigma\geq s_0\geq \rho^*$, we have that for any $\tau\in[\sigma,s]$, $\tau\leq s\leq 2\tau$.

   (i) Estimate of $\alpha_2(s)$.
    \begin{align}\alpha_2(s)
    =&\sigma^2s^{-2}q_2(\sigma)+(k_2,(\chi_1(\cdot,s)-\chi_1(\cdot,\sigma))\sigma^2s^{-2}q(\sigma))\nonumber\\
    &+(k_2,\chi_1(\cdot,s)(K(s,\sigma)-\sigma^2s^{-2})q(\sigma)).\end{align}
    By (\ref{qsigma}), we have
    $|\sigma^2s^{-2}q_2(\sigma)|\leq A''(\log \sigma)s^{-2} $, and \begin{align}|(k_2, (\chi_1(\cdot,s)-\chi_1(\cdot,\sigma))\sigma^2s^{-2}q(\sigma))|
\leq CA'(s-\sigma)s^{-3},\end{align}
for $\sigma\geq s_0\geq s_1^1(A',A'',A''',K_0,\rho^*)$.

We write $(K_2, \chi_1(\cdot,s)(K(s,\sigma)-\sigma^2s^{-2})q(\sigma))$ as $\sum_{r=0}^2 b_r+b_{-}+b_e$,
where $b_r=(k_2, \chi_1(\cdot,s)(K(s,\sigma)-\sigma^2s^{-2})h_r)q_r(\sigma)$, $b_{-}=(k_2, \chi_1(\cdot,s)(K(s,\sigma)-\sigma^2s^{-2})q_{-}(\sigma))$,
$b_{e}=(k_2, \chi_1(\cdot,s)(K(s,\sigma)-\sigma^2s^{-2})q_{e}(\sigma))$.

For $r=0,\ 1$, we use Lemma \ref{B-K1}, Corollary \ref{corollary 3.1} and (\ref{qsigma}), and the fact
that $e^{(s-\sigma)\mathcal{L}}h_r=e^{(1-r/2)(s-\sigma)}h_r$ and $(k_2,h_r)=0$, and derive
that when $s_0\geq s_1^2(A',K_0,\rho^*)$,
\begin{align}|b_r|=&|(k_2, \chi_1(\cdot,s)(K(s,\sigma)-e^{(s-\sigma)\mathcal{L}})h_r)q_r(\sigma)+(k_2,\chi_1(\cdot,s)(e^{(s-\sigma)\mathcal{L}}
-\sigma^2s^{-2})h_r)q_r(\sigma)|\nonumber\\
\leq  &CA'e^{(s-\sigma)}(s-\sigma)\sigma^{-3}. \label{br}\end{align}

Indeed,
by Lemma \ref{B-K1}, \begin{align}&|(K(s,\sigma)-e^{(s-\sigma)\mathcal{L}})h_r|\nonumber\\
\leq& \int |(K(s,\sigma,z,x)-e^{(s-\sigma)\mathcal{L}}(z,x))h_r(x)|dx\nonumber\\
\leq &\int e^{(s-\sigma)\mathcal{L}}(z,x)\Big( \frac{C}{s}(s-\sigma)\sum_{m,n\geq 0,m+n\leq2}|z|^m|x|^n\nonumber\\&+C(s-\sigma)(1+s-\sigma)s^{-2}(1+|z|+|x|)^4\Big)|h_r(x)|dx.\end{align}
Hence, by  Corollary \ref{corollary 3.1}  and (\ref{qsigma}), if $\sigma\geq s_0\geq s_1^3(\rho^*)$, then
\begin{align}&|(k_2, \chi_1(\cdot,s)(K(s,\sigma)-e^{(s-\sigma)\mathcal{L}})h_r)q_r(\sigma)|\nonumber\\
\leq &|\int k_2 \chi_1(z,s)(K(s,\sigma)-e^{(s-\sigma)\mathcal{L}})h_rq_r(\sigma)d\mu(z)|\nonumber\\
\leq &CA'(\frac{s-\sigma}{s^3}+\frac{(s-\sigma)(1+s-\sigma)}{s^4})e^{(s-\sigma)}
\leq CA'({s-\sigma})s^{-3}e^{(s-\sigma)}.\label{i4}\end{align}

On the other hand, just note that
\[
(e^{(s-\sigma)\mathcal{L}}
-\sigma^2s^{-2})h_r= (e^{(1-\frac r2)(s-\sigma)}- \sigma^2s^{-2})h_r.
\]
Note also that if $r\neq 2$, then $(k_2, h_r)=0$, hence
\[
(k_2, \chi_1(\cdot, s)h_r)= - (k_2, (1-\chi_1(\cdot,s))h_r),
\]
and
\[|(k_2, \chi_1(\cdot,s)
(e^{(s-\sigma)\mathcal{L}}
-\sigma^2s^{-2})h_r)q_r(\sigma)|=|(e^{(1-\frac r2)(s-\sigma)}- \sigma^2s^{-2})(k_2, (1-\chi_1(\cdot,s)h_r))||q_r(\sigma)|.
\]

 Remember first that we have
\[
s\in [\sigma, \sigma+\rho^*],
\]
and that we may assume that $\sigma \ge 1$.

\medskip

All we need to do, is to prove that
\begin{equation}\label{i}
\left|e^{(1-\frac r2)(s-\sigma)}
-\sigma^2s^{-2}\right|\le C(\rho^*)(s-\sigma).
\end{equation}
In order to do so, we write
\begin{equation}\nonumber
\left|e^{(1-\frac r2)(s-\sigma)}
  -\sigma^2s^{-2}\right|\le \left|e^{(1-\frac r2)(s-\sigma)}-1\right|+\left|1
-\sigma^2s^{-2}\right|\equiv I_1+I_2.
\end{equation}
Then, since for any $x\ge 0$, there is $y\in [0,x]$ such that
$|e^x-1|=e^yx\le xe^x$, we see that
\begin{equation}\label{i1}
I_1 \le e^{(1-\frac r2)(s-\sigma)}(s-\sigma) \le  e^{(1-\frac
  r2)\rho^*}(s-\sigma)\equiv C_1(\rho^*)(s-\sigma).
\end{equation}
As for $I_2$, we simply write
\begin{equation}\label{i2}
  I_2 = \frac{(s-\sigma)(s+\sigma)}{s^2}\le \frac{2s}{s^2}(s-\sigma)
  \le \frac 2\sigma(s-\sigma) \le 2(s-\sigma),
\end{equation}
where we have used the fact that $1\le \sigma \le s$.
Gathering \eqref{i1} and \eqref{i2}, we obtain \eqref{i}. Thus, when $s_0\geq s_1^4(A',\rho^*)$,
\begin{equation}|(k_2, \chi_1(\cdot,s)
(e^{(s-\sigma)\mathcal{L}}
-\sigma^2s^{-2})h_r)q_r(\sigma)|\leq Ce^{-cs}(s-\sigma),\nonumber
\end{equation}
from which and(\ref{i4}),  (\ref{br}) holds.

By Lemma \ref{B-K1} and (\ref{qsigma}),
$$|b_2|\leq CA'(s-\sigma)s^{-3},$$
if $\sigma\geq s_0\geq s_1^5(A',A'',K_0,\rho^*)$.

 We write
 $b_e=b_{e,1}+b_{e,2}+b_{e,3}$, with $b_{e,1}=(k_2, \chi_1(\cdot,s)(K(s,\sigma)-e^{(s-\sigma)\mathcal{L}})q_e(\sigma))
$,  $b_{e,2}=(k_2, \chi_1(\cdot,s)\int_0^{s-\sigma} d\tau \mathcal{L}(e^{\tau\mathcal{L}})q_e(\sigma))
$, $b_{e,3}=(k_2, \chi_1(\cdot,s)(1-\sigma^2s^{-2})q_e(\sigma))
$.
By (\ref{qsigma}),
$$|b_{e,3}|\leq C(s-\sigma)A's^{-3}
,$$ if $\sigma\geq s_0\geq s_1^6(A',A'',K_0,\rho^*)$.
Since $\mathcal{L}$ is self-adjoint,
\begin{align}|b_{e,2}|=&\Big|\int\frac{e^{-\frac{|z|^2}{4}}}{\sqrt{4\pi}}dz\mathcal{L}(k_2\chi_1(\cdot,s))(z)\nonumber\\&\cdot\int_0^{s-\sigma}d\tau\int dx\frac{e^\tau}{\sqrt{4\pi(1-e^{-\tau})}}
\exp\Big(-\frac{(ze^{-\tau/2}-x)^2}{4(1-e^{-\tau})}\Big)q_e(\sigma)\Big|\nonumber\\
\leq&\int\frac{e^{-\frac{|z|^2}{4}}}{\sqrt{4\pi}}dz\mathcal{L}(k_2\chi_1(\cdot,s))(z)\nonumber\\&\cdot\int_0^{s-\sigma}d\tau\int dx\frac{e^{\tau}}{\sqrt{4\pi(1-e^{-\tau})}}
\exp\Big(-\frac{(ze^{-\tau/2}-x)^2}{4(1-e^{-\tau})}\Big)A''\sigma^{-1/2},\nonumber\\
\leq &C(s-\sigma)A''s^{-1/2}e^{-Cs}\leq C(s-\sigma)A's^{-3},\end{align}
if   $\sigma\geq s_0\geq s_1^7(A',A'',K_0,\rho^*)$.

In the above inequality, we have used the fact that
$\mathcal{L}(k_2 \chi_1)$ is
zero for $|y|<K_0 \sqrt s$.
More precisely, since $\mathcal{L} k_2=0$ and $k_2'=k_1/2$, it is easy to see
that
\begin{align}
\mathcal{L}(k_2 \chi_1)
=\frac{z}{2}\chi_1'+\frac{z^2-2}{8}(\chi_1)''-\frac{1}{2}z\frac{z^2-2}{8}(\chi_1)'.
\end{align}
Since the support of the derivatives of $\chi_1$ is included in $[K_0
\sqrt s, 2K_0 \sqrt s]$, we derive that when $s_0\geq 1$,
\[
|\mathcal{L}(k_2 \chi_1) |\le C(K_0)(1+|z|^3)1_{\{|z|\ge K_0 \sqrt s\}}.
\]

Now, we estimate $b_{e,1}$.
By Lemma \ref{B-K1}, \begin{align}&|(K(s,\sigma)-e^{(s-\sigma)\mathcal{L}})q_e(\sigma)|\nonumber\\
\leq &\int |(K(s,\sigma,z,x)-e^{(s-\sigma)\mathcal{L}}(z,x))q_e(x,\sigma)|dx\nonumber\\
\leq &\int e^{(s-\sigma)\mathcal{L}}(z,x)\Big( \frac{C}{s}(s-\sigma)\sum_{m,n\geq 0,m+n\leq2}|z|^m|x|^n\nonumber\\&+C(s-\sigma)(1+s-\sigma)s^{-2}(1+|z|+|x|)^4\Big)|q_e(x,\sigma)|dx.\end{align}
Hence, by  Corollary \ref{corollary 3.1}  and (\ref{qsigma}), and the fact that $|z|\geq K_0\sqrt{\sigma}$, $q_e$ is not 0,  we have
\begin{align}&|b_{e,1}|\nonumber\\
\leq &\Big|\int k_2 \chi_1(z,s)(K(s,\sigma)-e^{(s-\sigma)\mathcal{L}})q_e(\sigma)d\mu(z)\Big|\nonumber\\
\leq &CA''\Big(\frac{s-\sigma}{s^{3/2}}+\frac{(s-\sigma)(1+s-\sigma)}{s^{5/2}}\Big)e^{(s-\sigma)}e^{-\frac{K_0^2}{8}\sigma}\nonumber\\
\leq &CA'({s-\sigma})s^{-3},\end{align}
if $\sigma\geq s_0\geq s_1^{8}(A',A'',K_0,\rho^*)$.

For
\begin{align}b_{-}=(k_2, \chi_1(\cdot,s)(K(s,\sigma)-\sigma^2s^{-2})q_{-}(\sigma)).\end{align}
 We write
 $b_{-}=b_{-,1}+b_{-,2}+b_{-,3}$, with $b_{-,1}=(k_2, \chi_1(\cdot,s)(K(s,\sigma)-e^{(s-\sigma)\mathcal{L}})q_-(\sigma))
$,  $b_{-,2}=(k_2, \chi_1(\cdot,s)\int_0^{s-\sigma} d\tau \mathcal{L}(e^{\tau\mathcal{L}})q_-(\sigma))
$, $b_{-,3}=(k_2, \chi_1(\cdot,s)(I-\sigma^2s^{-2})q_-(\sigma))
$.
Just note that
 $(k_2, q_-)=0$, we have
\[
(k_2, \chi_1(\cdot, s)q_-)= - (k_2, (1-\chi_1(\cdot,s))q_-).
\]
By (\ref{i2}) and (\ref{qsigma}),
\begin{align}|b_{-,3}|=&|(k_2, \chi_1(\cdot,s)(1-\sigma^2s^{-2})q_-(\sigma))|=|(1-\sigma^2s^{-2})||(k_2, (1-\chi_1(\cdot,s))q_-(\sigma))|
\nonumber\\\leq & C\frac{1}{s}(s-\sigma)A'''{\sigma^{-2}}e^{-Cs}\leq C(s-\sigma)A'''s^{-3}
,\end{align} if $\sigma\geq s_0\geq s_1^{9}(K_0,\rho^*)$.

$b_{-,2}$ can be treated similarly as $b_{e,2}$, it is bounded by
$C(s-\sigma)A'''s^{-2}e^{-Cs}\leq C(s-\sigma)A'''s^{-3},$
if   $\sigma\geq s_0\geq s_1^{10}(K_0,\rho^*)$.

By Lemma \ref{B-K1}, Corollary \ref{corollary 3.1}  and (\ref{qsigma}),
\begin{align}|b_{-,1}|=&|(k_2, \chi_1(\cdot,s)(K(s,\sigma)-e^{(s-\sigma)\mathcal{L}})q_-(\sigma)|\nonumber\\
\leq &|\int k_2, \chi_1(z,s)(K(s,\sigma)-e^{(s-\sigma)\mathcal{L}})q_-(\sigma)d\mu(z)|\nonumber\\
\leq & CA'''(\frac{s-\sigma}{s^3}+\frac{(s-\sigma)(1+s-\sigma)}{s^4})e^{(s-\sigma)}\nonumber\\
\leq &CA'''({s-\sigma})s^{-3}e^{(s-\sigma)},\end{align}
if $\sigma\geq s_0\geq s_1^{11}(\rho^*)$.

All these bounds yield
\begin{align}|\alpha_2(s)|\leq
  A''(\log \sigma)s^{-2} + C\max\{A',A'''\}(s-\sigma)e^{(s-\sigma)}s^{-3}, \end{align}
for $s_0\geq s_1^{12}(A',A'',A''',K_0,\rho^*)$.

(ii) Estimate of $\alpha_-(z,s)$.
\begin{align}\label{alpha(z,s)}\alpha_-(z,s)=&P_-(\chi_1(\cdot,s)K(s,\sigma)q(\sigma))\nonumber\\
=&P_-(\chi_1(\cdot,s)K(s,\sigma)q_-(\sigma))+\sum_{r=0}^2q_r(\sigma)P_-(\chi_1(\cdot,s)K(s,\sigma)h_r)\nonumber\\&+
P_-(\chi_1(\cdot,s)K(s,\sigma)q_e(\sigma)),\end{align}
where $P_-$ is the $L^2(\mathbb{R},d\mu)$ projector on the negative subspace of $\mathcal{L}$.

In order to bound the first term,
\begin{align}\label{Kq}
K(s,\sigma)q_-(\sigma)=\int dx e^{x^2/4}K(s,\sigma,z,x)f(x),
\end{align}
where $f(x)= e^{-x^2/4}q_-(x,\sigma)$.
From step 1, $e^{x^2/4}K(s,\sigma,z,x)=J(z,x)E(z,x)$
with
\begin{align}\label{N}
J(z,x)=[4\pi(1-e^{-(s-\sigma)}]^{-1/2}e^{s-\sigma}e^{x^2/4}e^{-\frac{(ze^{-(s-\sigma)/2}-x)^2}{4(1-e^{-(s-\sigma)})}}
\end{align}
and
\begin{align}
E(z,x)=\int d\mu_{zx}^{s-\sigma}(\omega)e^{\int_0^{s-\sigma}V(w(\tau),\sigma+\tau) d\tau}.
\end{align}

Let $f^0=f$, and for $m=0,1,2$, \[
  f^{(-m-1)}(z) = \int_{-\infty}^z   f^{(-m)}(x) dx,
\]
we have that for $m=0,1,2,3$ (see Lemma 6 in \cite{Bricmont}),
$$|f^{(-m)}(z)|\leq
CA^{'''}s^{-2}(1+|z|)^{3-m}e^{-\frac{z^2}{4}}.$$

Then, we can use the similar techniques used  in the proof of Lemma 3.13 in \cite{Zaag1} to get that
 if $s_0\geq s_1^{13}(A',A'',A''',K_0,\rho^*)$, then
$$|\alpha_-(z,s)|\leq C (A^{'''}s^{-2}e^{-(s-\sigma)/2}+A''e^{-(s-\sigma)^2} s^{-2})(1+|z|^3).$$

(iii) Estimate of $\alpha_e(z,s)$.

We write $$\alpha_e(z,s)=(1-\chi_1(z,s))K(s,\sigma)q(\sigma)=(1-\chi_1(z,s))K(s,\sigma)(q_b(\sigma)+q_e(\sigma)).$$
From (\ref{qsigma}) and Corollary \ref{corollary 3.1},
\begin{align}&|(1-\chi_1(z,s))K(s,\sigma)q_b(\sigma)|\leq Ce^{(s-\sigma)}A'''\sigma^{-1/2},\end{align}
if $\sigma\geq s_0\geq s_1^{14}(A',A'',A''',\rho^*)$.

Using (\ref{qsigma}) and the following estimate (see in \cite{Bricmont}),
$$\|K(s,\sigma)I_{|z|\geq K_0\sqrt{\sigma}}\|\leq Ce^{-(s-\sigma)/p},$$
we have
\begin{align}&|(1-\chi_1(z,s))K(s,\sigma)q_e(\sigma)|\nonumber\\=& \Big|(1-\chi_1(z,s))\int dx K(s,\sigma,z,x)q_e(x,\sigma)\Big|\nonumber\\
=& \Big|(1-\chi_1(z,s))\int dx K(s,\sigma,z,x)I_{|x|\geq K_0\sqrt{\sigma}}q_e(x,\sigma)\Big|\nonumber\\
\leq& C A''s^{-1/2}e^{-(s-\sigma)/p}.\end{align}
Hence, it holds that
$$\|\alpha_e(s)\|_{L^\infty}\leq Ce^{(s-\sigma)}A'''s^{-1/2}+C A''s^{-1/2}e^{-(s-\sigma)/p},$$
if $\sigma\geq s_0\geq s_1^{15}(A',A'',A''',\rho^*)$.

This completes the proof of Lemma \ref{BK}\ \ \ \#.\\ \\

\end{document}